\pdfoutput=1 
\documentclass{amsart}

\usepackage[pagewise]{lineno}
 
\usepackage{amsmath}
\usepackage{amssymb}
\usepackage{paralist}
\usepackage{graphics}
\usepackage{epsfig}
\usepackage{graphicx}  
\usepackage{epstopdf}
\usepackage{stmaryrd}

\usepackage[colorlinks=true]{hyperref}
\hypersetup{urlcolor=blue, citecolor=red}
\usepackage{hyperref}

\usepackage{rotating}
\usepackage{adjustbox}
\usepackage{array}      
\usepackage{multicol}
\usepackage[bottom]{footmisc}
\usepackage{float}
\usepackage{cancel}
\usepackage[usenames,dvipsnames]{pstricks}
\usepackage{multirow} 
\usepackage{color}
\usepackage{setspace}
\usepackage{colortbl}
\usepackage{caption}

\newtheorem{proposition}{Proposition}[section]

\numberwithin{equation}{section}

\newcommand{\R}{\mathbb{R}}

\newcommand{\f}[2]{\frac{\displaystyle{#1}}{\displaystyle{#2}}}
\newcommand{\n}{\newline}

\def\be{\begin{equation}}
\def\ee{\end{equation}}

\def\R{\mathbb{R}}

\def\f{\frac}

\def\ben{\begin{eqnarray}}
\def\een{\end{eqnarray}}

\def\ee{e^{\alpha z^*}}

\def\n{\newline}

\newcommand{\alain}[1]{{#1}}

\title[Bounded stochastic perturbations in population dynamics models]{\alain{A way to model stochastic perturbations in population dynamics models with bounded realizations}}
\author[T. Caraballo, R. Colucci, J. L\'opez-de-la-Cruz, A. Rapaport ]{}
				


\email{caraball@us.es}
\email{renatocolucci@hotmail.com}
\email{jlopez78@us.es}
\email{alain.rapaport@inra.fr}

\thanks{$^*$ Corresponding author: jlopez78@us.es}

\thanks{This work has been partially supported by FEDER and the Spanish  Ministerio de Econom\'{i}a y Competitividad project MTM2015-63723-P and the Consejer\'{i}a de Innovaci\'on, Ciencia y Empresa (Junta de Andaluc\'{i}a) under grant 2010/FQM314 and Proyecto de Excelencia P12-FQM-1492.}
	
\begin{document}

\maketitle

\centerline{\scshape Tom\'as Caraballo }
\medskip
{\footnotesize
 % please put the address of the second  and third author
 \centerline{Departamento de Ecuaciones Diferenciales y An\'alisis Num\'erico}
   \centerline{C/ Tarfia s/n, Facultad de Matem\'aticas, Universidad de Sevilla}
   \centerline{ 41080 Sevilla, Spain}
}
\medskip
 \centerline{\scshape Renato Colucci}
\medskip
{\footnotesize
 \centerline{Departament of Engineering,}
 \centerline{University Niccol\`{o} Cusano,}
   \centerline{ Via Don Carlo Gnocchi, 3
00166, Roma, Italy. }
}
\medskip
\centerline{\scshape Javier L\'opez-de-la-Cruz$^*$}
\medskip
{\footnotesize
\centerline{Departamento de Ecuaciones Diferenciales y An\'alisis Num\'erico}
   \centerline{C/ Tarfia s/n, Facultad de Matem\'aticas, Universidad de Sevilla}
   \centerline{ 41080 Sevilla, Spain}
}
\medskip
\centerline{\scshape and Alain Rapaport}
\medskip
{\footnotesize
 \centerline{MISTEA, Univ. Montpellier, Inra, Montpellier SupAgro}
   \centerline{2, place Pierre Viala}
   \centerline{34060 Montpellier, France}
}

\bigskip

\centerline{(Communicated by the associate editor name)}

\begin{abstract}
In this paper, we analyze the use of the Ornstein-Uhlenbeck process to model dynamical systems subjected to bounded noisy perturbations. In order to discuss the main characteristics of  this new approach we consider some basic models in population dynamics such as the logistic equations and competitive Lotka-Volterra systems. The key is the fact that these perturbations can be ensured to keep inside some interval that can be previously fixed, for instance, by practitioners, even though the resulting model does not generate a random dynamical system. However, one can still analyze the forwards asymptotic behavior of these random differential systems. Moreover, to illustrate the advantages of this type of modeling, we exhibit an example testing the theoretical results with real data, and consequently one can see this method as a realistic one, which can be very useful and helpful for scientists.
\end{abstract}
		
\vskip0.2cm

\section{Introduction}

%The mathematical modeling plays an important role when analyzing many problems in several branches of knowledge such as Biology, Chemistry, Physics and Sciences of life. In these domains, Applied Mathematics can provide useful techniques which allow us to obtain detailed information about the dynamics of the real models to be analyzed, either deterministic or stochastic. In this way, thanks to the theoretical analysis of the resulting mathematical models, it is possible to provide practitioners conditions on the parameters involved in such a models to ensure, for instance, the extinction of the populations of species that they need to analyze or, what is much more interesting from the mathematical point of view and important from the biological one, to provide conditions to ensure the persistence and the coexistence of the species in study.\n

In the last years, researchers of many areas in life sciences have been more and more interested in considering non-deterministic parameters in the mathematical models since it allows them to set up models which are much more realistic in the sense that they fit much better the ones in the real life. However, there are many different ways of introducing random or stochastic disturbances in deterministic models. Therefore, there are many details to be taken into account before starting, for instance

\begin{itemize}
\item {\bf Which kind of stochastic/random perturbation can we introduce?} There are many different stochastic processes and we need to make a decision about the one that we will use in order to set up our stochastic/random model. 
\item {\bf How can we do it?} Once decided the stochastic process to perturb the system, we need to think where and how we can introduce the disturbances.  
\item {\bf Is it realistic from the biological point of view?} After thinking about the previous questions, we should think whether our resulting stochastic/random model is realistic since our goal consists of obtaining models which reflect the reality as much better as possible.
\item {\bf Is it tractable from the mathematical point of view?} Of course, we need our models to be as much realistic as possible but we also need to have some tractability in order to make \alain{mathematical analysis and effective computation}, since no work could be made otherwise. Due to this facts, we need to find some suitable balance which make our work both tractable and realistic and, then, original and interesting.
\end{itemize} 

In the following part of this introductory section, we present a way of modeling random fluctuations in mathematical models of life sciences which will allow us to obtain realistic models in a simple way and, moreover, easy to treat from the mathematical point of view. In addition, we will test this way of modeling with real data in an specific model of a bioreactor to observe how it fits the perturbations that happen in the real experiment which represents a useful advance in the knowledge of the real problem.\n

The most common stochastic process that is considered when modeling disturbances in real life is the well-known standard Wiener process, see for instance \cite{cc,Graham2014} where the authors study random and stochastic modeling for a SIR model, \cite{cch,cch2,qx} \alain{where stochatic prey-predator Lotka-Volterra systems are analyzed} or \cite{CGL,CGLii,corrigendumchapter,imhof-walcher} \alain{where different ways of modeling stochasticity in the chemostat model are investigated}. Nevertheless, this stochastic process has the property of having continuous but not bounded variation paths, which does not suit to the idea of modeling real situations since, in most of cases, the real life is subjected to fluctuations which are known to be bounded.\n

It is worth mentioning that, in some cases analyzed in the previous literature (see \cite{CGL,corrigendumchapter}), the use of a standard Wiener process to perturb some parameters in deterministic systems to model random perturbations of different nature can lead to a non-realistic model, for instance, the positiveness of solutions are not necessarily preserved as a consequence of the arbitrary large values and the huge fluctuations of this Wiener process. However, in some cases one can modify the way in which the deterministic model can be replaced by a stochastic one, still using standard Wiener process, in such a way that the positiveness of solutions is preserved too (see \cite{CGLii,imhof-walcher}). From the mathematical point of view this way of modeling makes perfect sense although one may not be so sure from the biological one.\n

Henceforth, in this paper we consider a noise which remains bounded for every time in some interval to be previously fixed, for instance, by practitioners and, in addition, allows us to make calculations in a very simple way. Apart from that, we will present some examples to illustrate the advantages of using this noise when modeling random fluctuations, for example, we will be able to obtain realistic mathematical models due to the boundedness of the stochastic process to be used, we will be also able to prove the existence of absorbing and attracting sets which will not depend on the realization of the noise and we will be able to ensure the persistence and coexistence of the population of the species under some conditions on the parameters involved in the models.\n

In addition, every result will be proved forwards in time, differently to many other contributions in the literature where the authors need to set up a much more complicated framework based on the concepts of pullback attractor and pullback convergence within the theory of random dynamical systems. However, despite of these pullback concepts have proven very useful for the development of the theory of random dynamical systems, in some cases from applications it might not provide meaningful information about the forwards behavior of the system. We would also like to remark that, in this work, our resulting random systems will not generate a random dynamical system but this fact will not be a problem since it is possible to investigate the long-time behavior of the random system for every fixed event, which is a relevant improvement since it helps us to prove the forwards convergence that we use in our work.\n

The manuscript is organized as follows. In Section \ref{secou} we introduce a suitable Ornstein-Uhlenbeck (O-U) process depending on some parameters whose effects on the dynamics of the stochastic process are explained as well as some essential properties which will be very useful when making calculations in the following. In Sections \ref{s2} and \ref{s3} we present an example of a logistic differential equation with random disturbances in the environment and the growth rate by means of the O-U process. Therefore, in Section \ref{s4} we present an example concerning the parameter estimation in the logistic equation affected by the O-U process by setting up an {\it observer} which will consists of another differential system which will give information about the behavior of the state variables. Then, in Section \ref{s5} we present a random competitive Lotka-Volterra system where the growth rates of the species are affected by noise by means of the O-U process. In Section \ref{s6} we recall the advantages of using the O-U process when modeling reality in the case of the chemostat model. Finally, in Section \ref{s7} we present some comments and conclusions.

\section{The Ornstein-Uhlenbeck process.}\label{secou}

The key in our current work consists of perturbing the deterministic models by means of a suitable O-U process defined as the following random variable
\begin{equation}
z_{\beta,\gamma}^*(\theta _t\omega )=-\beta\gamma\int\limits_{-\infty }^0e^{\beta s}\theta _t\omega (s)ds,\quad t\in \ensuremath{\mathbb{R}},\,\, \omega \in \Omega,\,\, \beta,\,\,\gamma>0,\label{OU}
\end{equation}
\noindent where $\omega$ denotes a standard Wiener process in a certain probability space $(\Omega,\mathcal{F},\mathbb{P})$, $\beta$ and $\gamma$ are positive parameters which will be explained in more detail below and $\theta_t$ denotes the usual Wiener shift flow given by 
\begin{equation*}
\theta_t \omega(\cdot) = \omega(\cdot + t) - \omega(t),\quad t\in \ensuremath{\mathbb{R}}.\label{ws}
\end{equation*}

We note that the O-U process \eqref{OU} can be obtained as the stationary solution of the Langevin equation
\begin{equation}\label{eq:L}
dz +\beta z dt = \gamma d\omega.
\end{equation}

We would like to highlight that the O-U process is frequently used to transform stochastic models affected by the standard Wiener process into random ones (see \cite{CGL,CGLii,corrigendumchapter,cc}), which are much more tractable from the mathematical point of view, but both parameters $\beta$ and $\gamma$ are not taken into account or do not play any relevant role. Nevertheless, in the framework that we introduce in this paper we use directly this suitable O-U process depending on the parameters previously mentioned since they will be the key of the advantages provided by this way of modeling, as we will show in the rest of this work.\n

Due to the importance of those parameters, we introduce them now in a more detailed way and we show the effects that they have in the dynamics of the realizations of the O-U process.\n

The O-U process given by \eqref{OU} is a stationary mean-reverting Gaussian stochastic process where $\beta>0$ is a {\it mean reversion constant} that represents the strength with which the process is attracted by the mean or, in other words, how {\it strongly} our system reacts under some perturbation, and $\gamma>0$ is a {\it volatility constant} which represents the variation or the size of the noise independently of the amount of the noise $\alpha>0$. In fact, the O-U process can describe the position of some particle by taking into account the friction, which is the main difference with the standard Wiener process and makes our perturbations to be a better approach to the real ones than the ones obtained when using simply the standard Wiener process. In addition, the O-U process could be understood as a generalization of the standard Wiener process as well in the sense that it would correspond to take $\beta=0$ and $\gamma=1$ in \eqref{OU}. In fact, the O-U also provides a link between the standard Wiener process and no noise at all, as we will see later.\n

Now, we would like to illustrate the relevant effects caused by both parameters $\beta$ and $\gamma$ on the evolution of realization of the O-U process.\n

\begin{figure}[H]
\begin{center}
\includegraphics[width=1.0\textwidth]{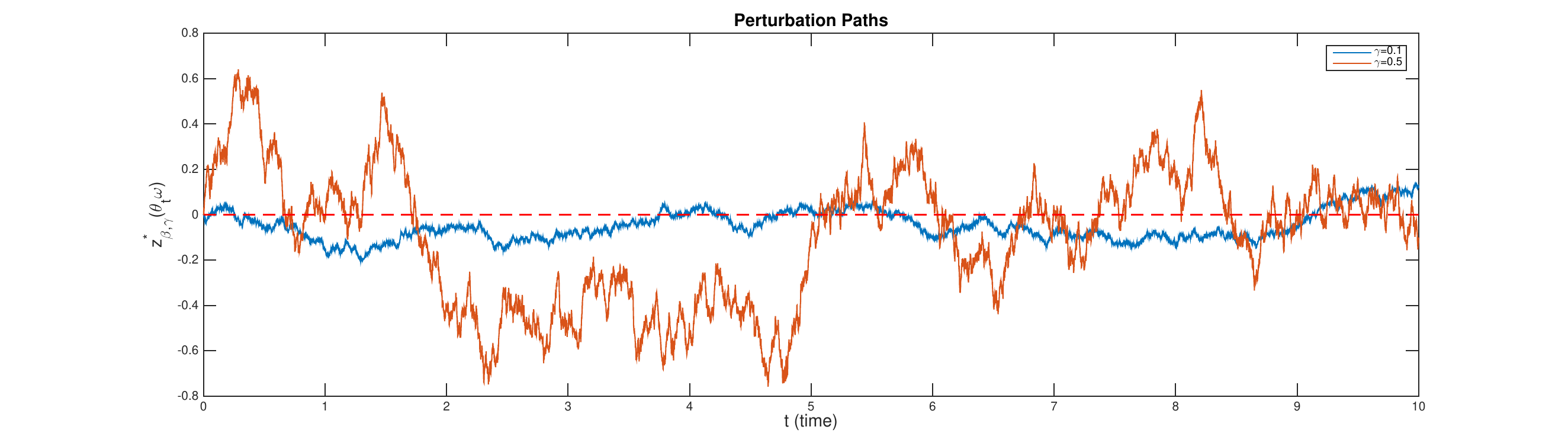}
\caption{Effects of the mean reverting constant on the O-U process}
\label{CGLR-p1}
\end{center}
\end{figure}

{\bf Fixed $\beta>0$.} Then, the volatility of the process increases when considering larger values of $\gamma$ and the evolution of the process is smoother when taking smaller values of $\gamma$, which sounds reasonable due to the fact that $\gamma$ decides the amount of noise introduced to $dz$, the term which measures the variation of the process. Henceforth, the process will be subject to suffer much more disturbances when taking a larger value of $\gamma$. This behavior can be observed in Figure \ref{CGLR-p1}, where we simulate two realizations of the O-U process with $\beta=1$ and we consider $\gamma=0.1$ (blue) and $\gamma=0.5$ (orange).

\begin{figure}[H]
\begin{center}
\includegraphics[width=1.0\textwidth]{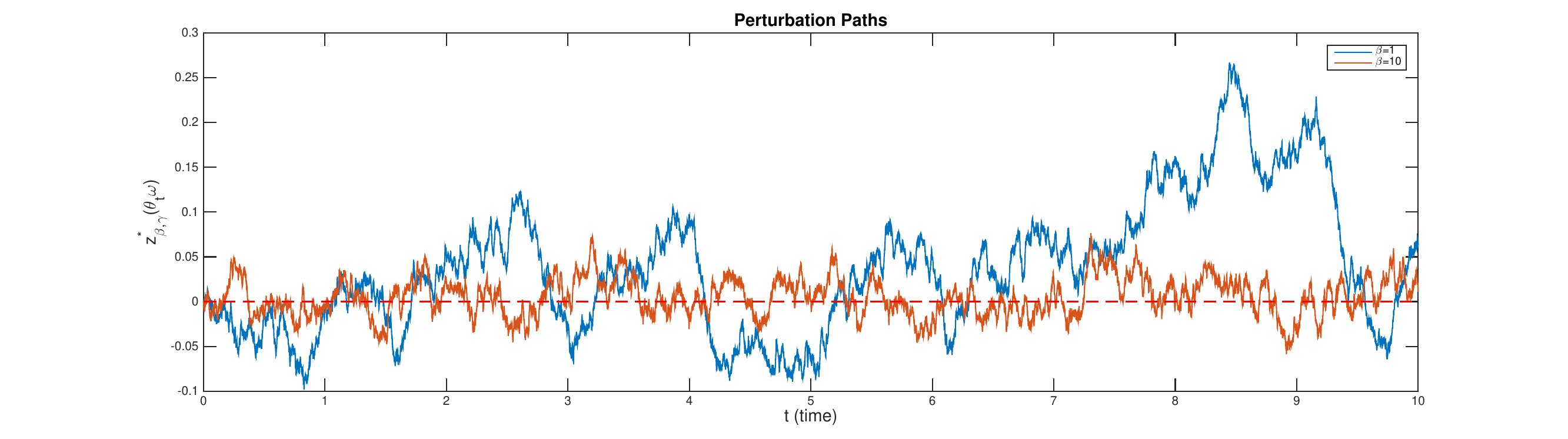}
\caption{Effects of the volatility constant on the O-U process}
\label{CGLR-p2}
\end{center}
\end{figure}

{\bf Fixed $\gamma>0$.} In this case the process tends to go further away from the mean value when considering smaller values of $\beta$ and the attraction of the mean value increases when taking larger values of $\beta$. This behavior seems logical since $\beta$ has a huge influence on the drift of the Langevin equation \eqref{OU}. We can observe this behavior in Figure \ref{CGLR-p2}, where we simulate two realizations of the O-U process with $\gamma=0.1$ and we take $\beta=1$ (blue) and $\beta=10$ (orange).

Once presented the O-U process and the effects that its parameters cause on \alain{the behavior of its realization}, we state now some essential properties that it satisfies which will be another important \alain{key point of the paper}.

\begin{proposition}\label{pOU} 
There exists a $\theta _t$-invariant set $\widetilde{\Omega }\in \mathcal{F}$  of  $\Omega$  of full $\ensuremath{\mathbb{P}}-$measure such that for $\omega \in \widetilde{\Omega }$ and $\beta,\gamma>0$,  we have
\begin{itemize}
\item [(i)] the random variable  $|z_{\beta,\gamma}^*(\omega )|$  is tempered.\n
\item [(ii)] the mapping 
\[
(t,\omega )\rightarrow z_{\beta,\gamma}^*(\theta _t\omega )=-\beta\gamma\int\limits_{-\infty
}^0e^{\beta s}\omega (t+s)\mathrm{d}s+\omega(t)
\]
is a stationary solution of \eqref{OU}
with continuous trajectories;\n
\item [(iii)] for any $\omega \in \tilde \Omega$ one has
\begin{eqnarray*}
\lim_{t\rightarrow \pm \infty }\frac{|z_{\beta,\gamma}^*(\theta _t\omega )|}%
t&=& 0;\\
\lim_{t\rightarrow \pm \infty }\frac 1t\int_0^tz_{\beta,\gamma}^*(\theta _s\omega
)ds&=&0;\\
 \lim_{t\rightarrow \pm \infty }\frac 1t\int_0^t |z_{\beta,\gamma}^*(\theta _s\omega
)| ds&=& \mathbb{E}[z_{\beta,\gamma}^*] < \infty;\n
\end{eqnarray*}
\item [(iv)] finally, for any $\omega\in\widetilde{\Omega}$,
$$\lim_{\beta\rightarrow\infty}z^*_{\beta,\gamma}(\theta_t\omega)=0, \quad\text{for all}\,\, t\in\R.$$
\end{itemize}
\end{proposition} 

The proof of Proposition \ref{pOU} is omitted here. We refer the readers to \cite{alazzawi2017} (Lemma 4.1) for the proof of the last statement and to \cite{arnold,caraballo1} for the proof of the other points.\n

To sum up the main idea of this framework, we will have to deal with a random differential system depending on the stationary solution of the Langevin equation \eqref{eq:L} as follows
\begin{equation}\label{eq:main}
\dot{x}=f(x, z_{\beta,\gamma}^*(\theta_t\omega)),
\end{equation}
\noindent where $z_{\beta,\gamma}^*(\theta_t\omega)$ is the stationary solution of the Langevin equation.\n

Henceforth, the idea consists of fixing \alain{an event} $\omega\in\Omega$ and then to solve \eqref{eq:L} in order to introduce its solution into the differential system \eqref{eq:main}. However, we remark that, for a general choice of the parameter $\beta$, it is not possible to ensure that the solution $z_{\beta,\gamma}^*(\theta_t\omega)$ is bounded for every time, which would be the most desirable property since it is well known that real models are \alain{expected to be subject} to bounded perturbations, in fact, unbounded disturbances has no sense when modeling \alain{real systems}, as we pointed out previously.\n

Nevertheless, thanks to the last property (iv) in Proposition \ref{pOU}, for every fixed event $\omega\in\Omega$, it is possible to find $\beta>0$ large enough such that the stationary solution of the Langevin equation \eqref{eq:L}, $z_{\beta,\gamma}^*(\theta_t\omega)$, remain inside any bounded interval previously fixed.\n

Therefore, for every fixed $\omega\in\Omega$ and for proper $\beta>0$, we compute $z_{\beta,\gamma}^*(\theta_t\omega)$ such that the realizations of the corresponding perturbed parameter remain inside a strictly positive interval, namely $(b_1,b_2)\subset\R$, which should be previously determined \alain{depending on the application or provided by practitioners}, where $b_2>b_1>0$. \alain{For population dynamics, this} property concerning the positiveness of the realizations of the perturbed parameter will be essential in this work, in fact, it will be the main key to be able to guarantee strictly positiveness of all state variables which could mean, for instance, to be able to ensure the persistence of the species involved in the corresponding model.\n

In most of situations one \alain{design a suitable} mathematical \alain{model with random features} for a biological problem. But the main concerns are: is the model realistic? is the model useful? To answer these questions one may try to \alain{identifty the parameters on sets of real data to test if the choice of stochastic process is able to reproduce satisfactorily real situations}.\n

\begin{figure}[H]
\begin{center}
\includegraphics[scale=0.25]{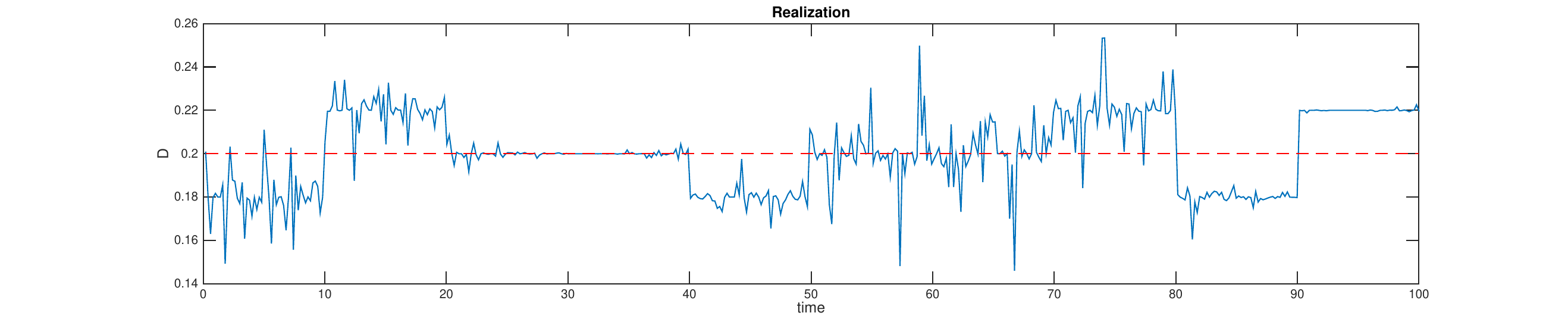}
\end{center}
\caption{Real data: input flow in a bioreactor}
\label{real}
\end{figure}

More precisely, in our case, the parameters $\beta$ and $\gamma$ in the O-U process can be estimated from the real data by using a simple mean square regression. As an example of how the process describes the fluctuations in the real models, in Figure \ref{real} we can see \alain{time history of} values that the input flow of a bioreactor has \alain{over a time horizon} of one-hundred hours. As it was expected, it seems to be mean-reverting and then it could be described as a realizations of the O-U process. When estimating the parameters by using a mean square regression, we obtain that the mean value is $\mu=0.2002$, the mean-reverting constant is $\beta=1.3230$ and the volatility constant is $\gamma=0.0287$.

\begin{figure}[H]
\begin{center}
\includegraphics[scale=0.25]{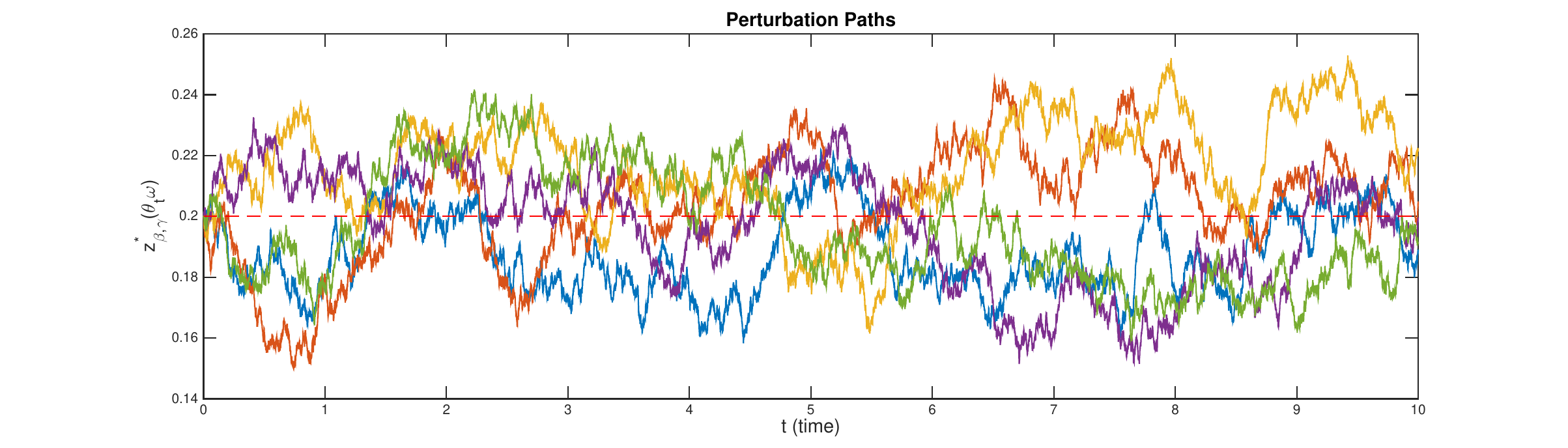}
\end{center}
\caption{Realizations of the O-U process generated with parameters from the real data with $\beta$, $\mu$, $\gamma$ as above}
\label{OUdata}
\end{figure}

Then, in Figure \ref{OUdata} we make some simulations of the O-U process with the previous values of the parameters obtained from Figure \ref{real}. As we can observe, every realization remains in the same interval $[0.14, 0.26]$.

Then, with biological applications in mind, we could consider perturbations of some parameter in a random system of the following kind
\begin{equation}\label{eq:main}
\dot{x}=f(x,z^*_{\beta,\gamma}(\theta_t\omega)).
\end{equation}

In the real world, one is often interested in analyzing bounded perturbations, thus one can wonder about a suitable way of modeling them. In our method for each $\omega$ fixed, we can find values of parameter $\beta_\omega$ such that the random perturbation is bounded in a prefixed interval. This implies that for each realization, we have that our system \eqref{eq:main} looks like
\begin{equation}
\dot{x}=f(x,z^*_{\beta,\gamma}(\theta_t\omega))=f(x,g(t))
\end{equation}
\noindent where $g(\cdot)$ is a continuous function for each $\omega$, and can be analyzed by the deterministic theory. The main fact is that in some cases the asymptotic behavior of this system may be independent of $\omega$, which then allows us to compare the deterministic and stochastic models and testing how realistic the latter can be.\n

In other words, the resulting random model is a non-autonomous system perturbed by means of a non-autonomous perturbation which is generated in a random way. As a consequence, in this context, the set of admissible random perturbation coincide with the set of continuous function satisfying properties (i)-(iv) stated in Proposition \ref{pOU}.\n

This kind of non-autonomous perturbations may be obtained by more general stochastic processes, then one could wonder the reason why we are just focusing on the O-U process instead of using another one. The answer is as clear as easy: \alain{Indeed,} this suitable O-U process given by \eqref{OU} provides us essential ergodic properties which allow us to make calculations when analyzing the resulting random systems, since we will have to deal with integral terms involving such a perturbation.\n

Notice that real models are subjected to fluctuations which are \alain{some what smooth w.r.t. time} and it would not make \alain{sense to have} perturbations which are \alain{changing too rapidly} from extreme values in short periods of time. This is \alain{one of the main reasons to consider } the O-U process by taking into account both parameters $\beta$ and $\gamma$ which allow us to \alain{have more flexibility to control} the noise in order to be able to represent the reality \alain{in a better way}.\n

One could also wonder the difference between this way of modeling with the one when considering some random function $a(\theta_t\omega)\in[a_1,a_2]$, which is bounded and continuous respect to the time, as in \cite{cc,caraballo5}. On the one hand, in both cases we obtain bounded perturbations, which is the most realistic from the point of view of the applications, as we explained previously. On the other hand, in both cases we are working with continuous functions respect to the time, which is also logical when trying to deal with differential systems as in the current case. Nevertheless, we remark that the random constant $a(\theta_t\omega)$ does not satisfy in general the properties (i)-(iv) stated in Proposition \ref{pOU}, as it happens with the O-U process, which are essential when making calculations and allow us to obtain \alain{characterizations of} absorbing and attracting sets for the solutions of the corresponding random systems as well. Moreover, we would have to assume $a(\theta_t\omega)$ to be continuous and bounded while the O-U process satisfies these properties by definition.\n

We would also like to remark that, in the classic random case, the continuous function $a(\theta_t\omega)$ is directly generated thanks to the dynamics of the set of events $\Omega$ whereas, in our new framework, every event $\omega\in\Omega$ is fixed and the continuous perturbation is obtained by solving the Langevin equation \eqref{eq:L}. Apart from that, the continuous function $a(\theta_t\omega)$ is an arbitrary continuous function with values in some positive bounded interval $[a_1,a_2]$ whereas, in our new framework, we are simply considering the realizations of the perturbed parameter which are realistic from the point of view of the applications.\n

\begin{center}
\begin{figure}[H]
\includegraphics[width=0.7\textwidth]{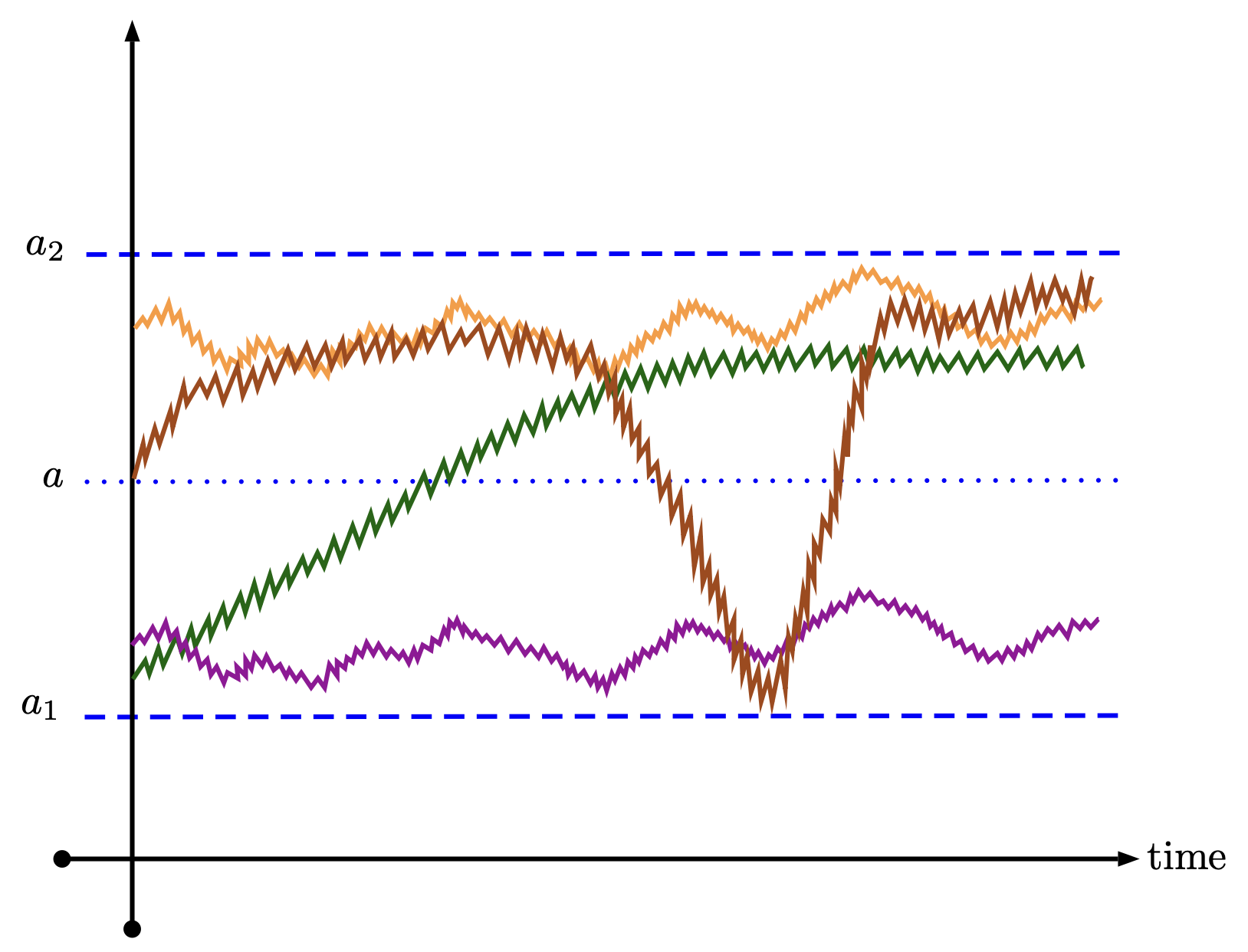}
\caption{Examples of non-realistic realizations of the perturbed parameter}
\label{norealperturbs}
\end{figure}
\end{center}

\alain{Let us underline that the realizations by means of suitable O-U process look quite different from the one generated by a random function $a(\theta_t\omega)$. One can observe time to time realizations of $a(\theta_t\omega)$ which are very unlikely to be observed as generating by the suitable O-U processes. Several kind of such realizations are depicted on Figure \ref{norealperturbs}:

\begin{itemize} 
\item approaching the boundaries of the interval $[a_1,a_2]$ and stay close to it for a while, as the green one
\item staying close from one of the boundaries, as the orange or violet ones
\item switching very rapidly between two values close from the boundaries of the interval, as the brown one.
\end{itemize}

The realizations generated by the suitable O-U looks more realistic in the sense that its is similar to an agitated particule with a recall force to the mean value.}

In both cases, when considering a \alain{classical random} \alain{function} $a(\theta_t\omega)$ and also when considering a perturbation by means of the O-U process $a+\alpha z^*_{\beta,\gamma}(\theta_t\omega)$, the perturbations are bounded, then one can expect to find bounds for the solutions of the system and, therefore, to be able to provide some conditions under which the persistence of the \alain{populations} involved in the model can be ensured. Nevertheless, there are important differences between both cases the classic and the new one, for instance, the natural context in the classic random case is to study the pullback convergence whereas, in the new random case involving the new suitable O-U process, the solutions may not generate a random dynamical system, since $\beta$ in fact depends on $\omega$. \alain{However, it} does not present any inconvenient since we can analyze the random system for every fixed $\omega\in\Omega$, as we explained before.\n

In addition, we can prove every mathematical result to hold forwards in time, which is much more realistic than the pullback convergence obtained in the classic case. This improvement concerning the forwards convergence is also very related with the ergodic properties stated Proposition \ref{pOU} (iii) which are proved to hold also forwards in time.\n

In order to illustrate \alain{the modeling approach we propose}, we are going to consider in this paper two \alain{well known models in the ecology literature: the logistic and the competition Lotka-Volterra model, that we revisit here introducing noise.}. \alain{In addition}, we will consider an observer \alain{dynamics} to \alain{show the parameter estimation despite noise} in the logistic equation affected by the O-U process, a problem that has not been treated up to now in the literature even though it is object of very interest.

\section{Environmental perturbations in the logistic model.}
\label{s2}

In this section we consider the classical logistic equation given by
\begin{eqnarray}
\label{syslogistic}
\f{dx}{dt}&=&x\left(a-x\right),
\end{eqnarray}
\noindent where $x=x(t)$ denotes the number of individual of some population of a certain species and $a$ is the carrying capacity.\n

As it is very well-known, the carrying capacity of the population whose dynamics are modeled by the logistic equation can be affected by many \alain{external factors} present on the environment as the climate or the temperature, to name a few. Then it has sense to consider random perturbations on the carrying capacity such that we have the following random equation
\begin{eqnarray}
\f{dx}{dt}&=&x\left(a+\alpha z^*_{\beta,\gamma}(\theta_t\omega)-x\right),\label{logisticrenato}
\end{eqnarray}
\noindent where $z^*_{\beta,\gamma}(\theta_t\omega)$ denotes the Ornstein-Uhlenbeck process that we introduced previously and $\alpha>0$ is the amount of noise.\n

The solution of the random logistic equation \eqref{logisticrenato} exists and its explicit expression is given by
\begin{equation}
x(t;0,\omega,x_0)=\f{x_0e^{\int_0^t a+\alpha z^*_{\beta,\gamma}(\theta_s\omega)ds}}{1+x_0\int_0^te^{\int_0^s a+z^*_{\beta,\gamma}(\theta_\tau\omega)d\tau}ds},\label{solutionrenato}
\end{equation}
for any initial value $x_0\geq 0$, any $\omega\in\Omega$ and for all $t\geq 0$.\n

In addition, thanks to a suitable choice of the parameter $\beta$ in the O-U process presented in the introduction of the paper, we know that , $a+\alpha z^*_{\beta,\gamma}(\theta_t\omega)\in[\underline{a},\bar{a}]$ for every $t\in\R$, \alain{where $\overline a > \underline a$ are positive values}. Then from \eqref{logisticrenato} we can obtain the following differential inequalities
\begin{equation}
x\left(\underline{a}-x\right)\leq \f{dx}{dt}\leq x\left(\bar{a}-x\right),\label{inrenato}
\end{equation}
\noindent whence we can deduce that, as soon as we consider an initial value of the species $x_0<\underline{a}$, then the dynamics of the population is increasing till it reaches the curve $a+\alpha z^*_{\beta,\gamma}(\theta_t\omega)$ which remains inside the positive interval $[\underline{a},\bar{a}]$.\n

Henceforth, from \eqref{inrenato} it can be deduced that, for every $\varepsilon>0$, any $\omega\in\Omega$ and any initial value $x_0<\underline{a}$, there exists some time $T(\varepsilon,\omega)>0$ such that
\begin{equation}
\underline{a}-\varepsilon\leq x(t;0,\omega,x_0)\leq \bar{a}+\varepsilon,
\end{equation}
\noindent for every $t\geq T(\varepsilon,\omega)$.\n

From the previous analysis we obtain that, for any $\varepsilon>0$, $B_\varepsilon=[\underline{a}-\varepsilon,\bar{a}+\varepsilon]$ is a deterministic absorbing set for the solutions of \eqref{logisticrenato}.\n

Therefore, $B_0=[\underline{a},\bar{a}]$ is a positive attracting set for the solutions of \eqref{logisticrenato}, i.e.,
\begin{equation}
\lim_{t\rightarrow+\infty}\sup_{x_0\in(0,\underline{a})}\inf_{b_0\in B_0}|x(t,0,\omega,x_0)-b_0|=0.
\end{equation}

Now, we present some numerical simulations to support the results previously provided and the advantages of using the suitable O-U process presented here when modeling realistic problems. From now on, the blue dashed lines represent the solutions of the deterministic models and the rest are different realizations of the random ones.\n

\begin{figure}[H]
\begin{center}
\includegraphics[scale=0.25]{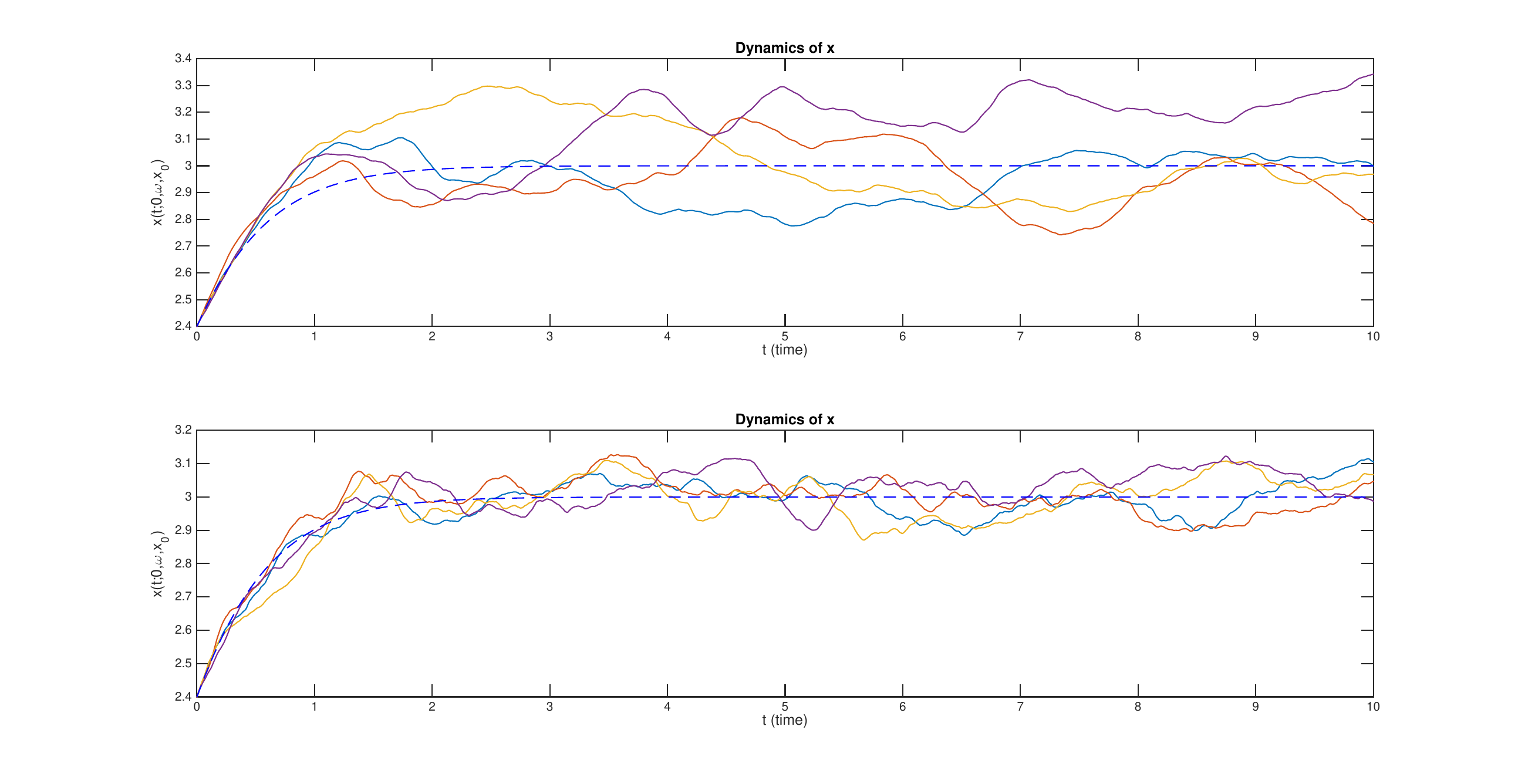}
\end{center}
\caption{Realizations of the solution of the random logistic equation with perturbed carrying capacity for $x_0=2.4$.  $\alpha=2$, $\beta=1$, $\gamma=0.1$(top) and $\alpha=2.2$, $\beta=10$, $\gamma=0.2$ (bottom)}
\label{renato1}
\end{figure}

In Figure \ref{renato1} we can see two panels representing several realizations of the solution of the random logistic equation \eqref{logisticrenato} for the initial value $x_0=2.4$, the nominal carrying capacity is $a=3$, the amount of noise is $\alpha=2$ (top) and $\alpha=2.2$ (bottom), the mean reverting constant is $\beta=1$ (top) and $\beta=10$ (bottom) and the volatility constant is $\gamma=0.1$ (top) and $\gamma=0.2$ (bottom).

\begin{figure}[H]
\begin{center}
\includegraphics[scale=0.25]{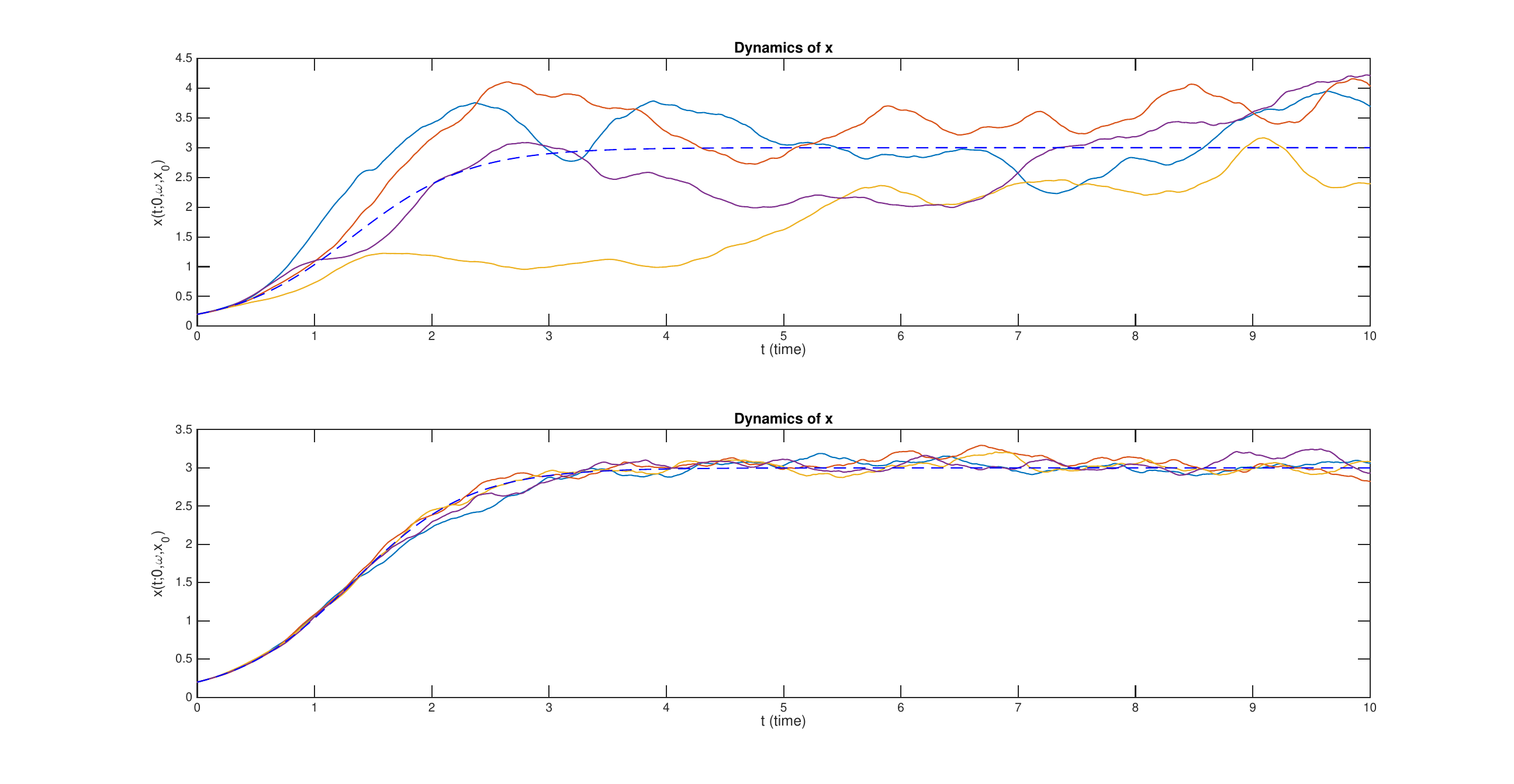}
\end{center}
\caption{Realizations of the solution of the random logistic equation with perturbed carrying capacity for $x_0=0.2$. $\alpha=2$, $\beta=1$, $\gamma=0.4$(top) and $\alpha=2.2$, $\beta=10$, $\gamma=0.4$ (bottom)}
\label{renato4}
\end{figure}

In Figure \ref{renato4} we display two panels representing several realizations of the solution of the random logistic equation \eqref{logisticrenato} for the initial value $x_0=0.2$, the nominal carrying capacity is $a=3$, the amount of noise is $\alpha=2$ (top) and $\alpha=2.2$ (bottom), the mean reverting constant is $\beta=1$ (top) and $\beta=10$ (bottom) and the volatility constant is $\gamma=0.4$. \alain{Compared to} Figure \ref{renato1}, now we increase the volatility constant which is significant for small values of the mean reverting constant (as it can be seen in the figure of the top) but the noise can be reduced if we increase the mean reverting constant even though the volatility constant is not decreased.

We can observe that all the solutions of the random equation \eqref{logisticrenato} are fluctuating around the equilibrium of the deterministic case $x=3$ and these fluctuations remain inside a strictly positive bounded interval which is smaller when taking larges values of $\beta$ and (or) smaller values of $\gamma$. Thus, the theoretical results and the advantages of the O-U process are \alain{demonstrated on this example}.\n

\begin{figure}[H]
\begin{center}
\includegraphics[scale=0.25]{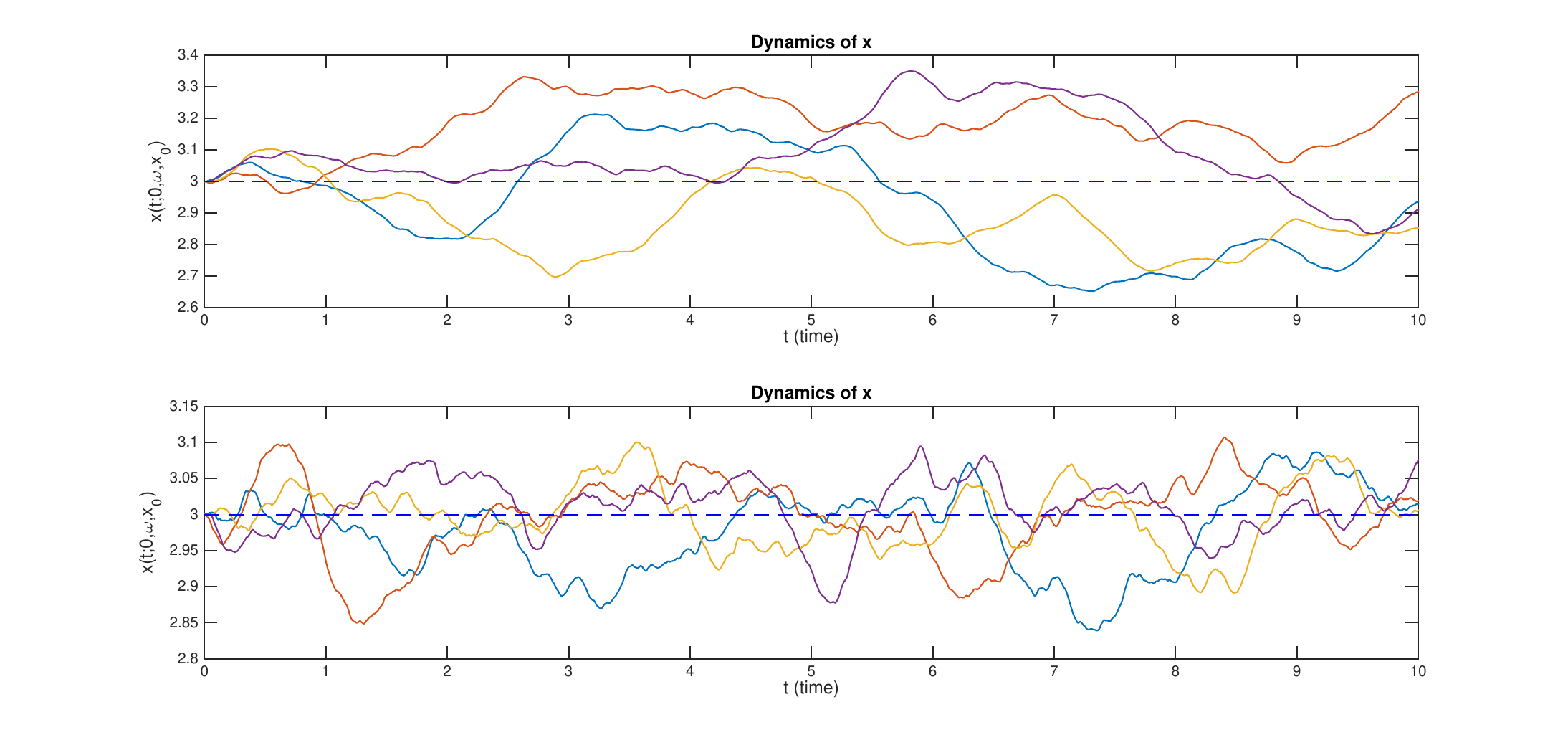}
\end{center}
\caption{Realizations of the solution of the random logistic equation with perturbed carrying capacity for $x_0=3$. $\alpha=2$, $\beta=1$, $\gamma=0.1$(top) and $\alpha=2.2$, $\beta=10$, $\gamma=0.2$ (bottom)}
\label{renato2}
\end{figure}

In Figure \ref{renato2} we present the behavior of several realizations of the solution of the random logistic equation \eqref{logisticrenato} for the initial value $x_0=3$, the nominal carrying capacity is $a=3$, the amount of noise is $\alpha=2$ (top) and $\alpha=2.2$ (bottom), the mean reverting constant is $\beta=1$ (top) and $\beta=10$ (bottom) and the volatility constant is $\gamma=0.1$ (top) and $\gamma=0.2$ (bottom).

We can observe in this case a similar behavior to the previous one. However \alain{there are significant} differences when comparing the behavior of the random equation \eqref{logisticrenato} and the deterministic one for the initial condition $x_0=3$. In the deterministic case, $x=3$ is a stable equilibrium and the solution is constant as can be observed in the blue dashed lines. \alain{On the random case, one can observe} some fluctuations around $x=3$ which remain inside a strictly positive interval, this later one being deterministic, as proved in the theoretical results. In addition, this interval can be chosen by tuning the parameters of the O-U process as we previously explained.

\section{Perturbations on the growth rate in the logistic equation.}\label{s3} 

Now, we consider the logistic equation that we rewrite in the following form
\begin{eqnarray}
\f{dx}{dt}&=&rx\left(1-\f{x}{c}\right),
\end{eqnarray}
\noindent where $x=x(t)$ denotes the number of population of some species, $r$ denotes the specific growth rate of the species and $c$ is the carrying capacity of the medium assumed to be constant, both positive.\n

In this case we are interested in introducing a noise in the reproduction rate by using the O-U process. As a result, we have the following random logistic model
\begin{eqnarray}
\f{dx}{dt}&=&(r+\alpha z^*_{\beta,\gamma}(\theta_t\omega))x\left(1-\f{x}{c}\right),\label{logistical}
\end{eqnarray}
\noindent where $z^*_{\beta,\gamma}(\theta_t\omega)$ denotes again the O-U process and $\alpha>0$ is the amount of noise. We observe that $x=c$ is still an equilibrium for the equation.\n

As made in the previous case, the solution of equation \eqref{logistical} exists and its explicit expression is given by
\begin{equation}
x(t;0,\omega,x_0)=\f{x_0}{e^{-\int_0^t r+\alpha z^*_{\beta,\gamma}(\theta_s\omega)ds}\left(1+\f{x_0}{c}\right)+\f{x_0}{c}}\label{solutional}
\end{equation}
\noindent for every $x_0\geq 0$, any $\omega\in\Omega$ and $t\geq 0$, whence we observe the property
\begin{eqnarray}
\nonumber
\int_0^t r+\alpha z^*_{\beta,\gamma}(\theta_s\omega)ds &=& rt+\int_0^t z^*_{\beta,\gamma}(\theta_s\omega) ds
\\[1.3ex]
\nonumber
&=&t\left(r+\f{1}{t}\int_0^t z^*_{\beta,\gamma}(\theta_t\omega)ds\right)
\end{eqnarray}

Thus, thanks to the ergodic properties in Theorem \ref{pOU}, we obtain that the dynamics of the population converges to the carrying capacity as in the deterministic case or, in other words, we have that for every $\varepsilon>0$, any initial value $0<x_0<c$ and $\omega\in\Omega$ , there exists some time $T(\varepsilon,\omega)>0$ such that
\begin{equation}
c-\varepsilon<x(t;0,\omega,x_0)<c
\end{equation}
\noindent for all $t\geq T(\varepsilon,\omega)$.\n

Therefore, we have that $B_\varepsilon=[c-\varepsilon,c]$, for any $\varepsilon>0$, is a deterministic absorbing set for the solutions of \eqref{logistical} whence we have that $B_0=\{c\}$ defines a positive deterministic attracting set for the solutions of \eqref{logistical}. As a consequence, every realization of the solution of \eqref{logistical} converge to the carrying capacity $c$ as long as the initial value $x_0>0$, as in the deterministic case. This is not surprising since, in this second logistic equation, the carrying capacity is still a stable equilibrium even though we are treating a random case.\n

We would also like to remark that the above calculations are independent
of the choice of the parameter $\beta$.\n

\begin{figure}[H]
\begin{center}
\includegraphics[scale=0.25]{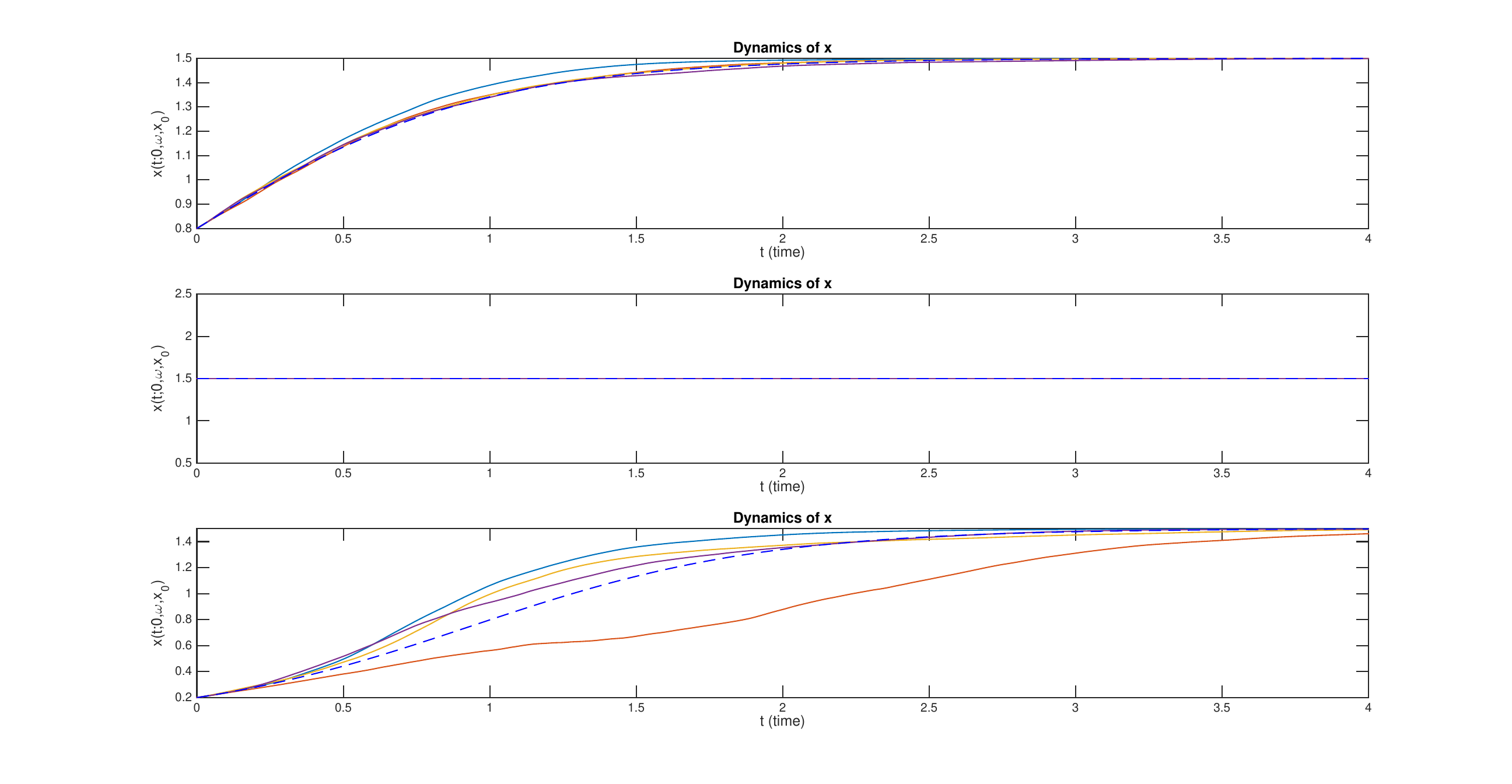}
\end{center}
\caption{Realizations of the solution of the random logistic equation with perturbed growth rate for $x_0=0.8$ (top), $x_0=1.5$ (medium) and $x_0=0.2$ (bottom)}
\label{alain1}
\end{figure}

Now we perform some numerical simulations to support the results previously stated. In Figure \ref{alain1} we show the behavior of several realizations of the solution of the random logistic equation \eqref{logistical} for the initial values $x_0=0.8$ (top), $x_0=1.5$ (medium) and $x_0=0.2$ (bottom), the growth rate is $r=2$, the carrying capacity is $c=1.5$, the amount of noise is $\alpha=2$, the mean reverting constant is $\beta=1$ and the volatility constant is $\gamma=0.4$.

We can observe, differently to the example analyzed in Section \ref{s2}, that in this case the realizations of the solution of the random equation \eqref{logistical} have fluctuations when the population is increasing but these disturbances are not present when the population is close to the carrying capacity, in fact, $x=c$ is an equilibrium of the equation \eqref{logistical} as in the deterministic case. In the second plot we can in fact see that the random solution is constant for the initial condition $x_0=1.5$.\n

\begin{figure}[H]
\begin{center}
\includegraphics[scale=0.25]{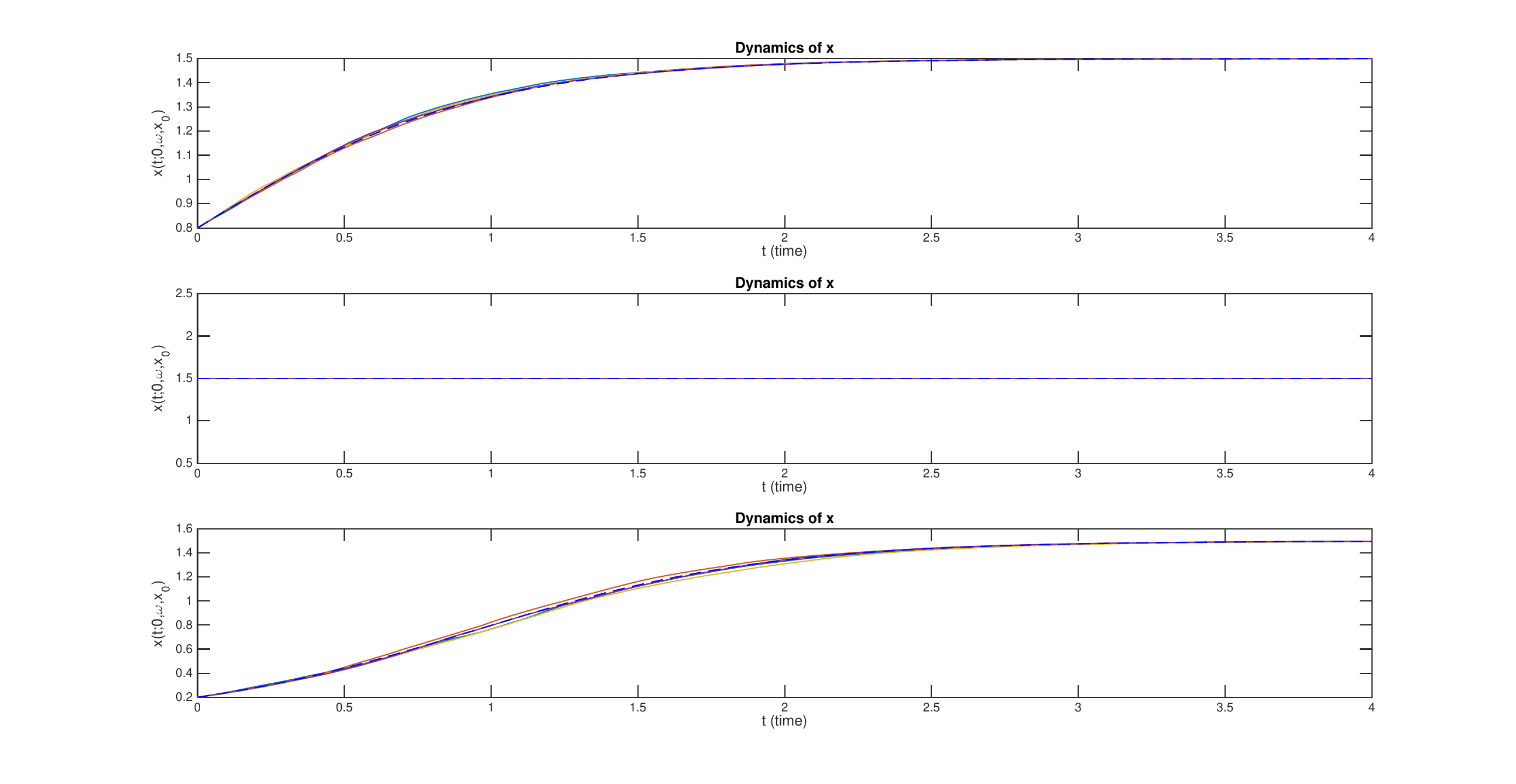}
\end{center}
\caption{Realizations of the solution of the random logistic equation with perturbed growth rate for $x_0=0.8$ (top), $x_0=1.5$ (medium) and $x_0=0.2$ (bottom)}
\label{alain2}
\end{figure}

In Figure \ref{alain2} we can see two panels representing several realizations of the solution of the random logistic equation \eqref{logistical} where the initial values $x_0=0.8$ (top), $x_0=1.5$ (medium) and $x_0=0.2$ (bottom), the growth rate is $r=2$, the carrying capacity is $c=1.5$, the amount of noise is $\alpha=2$, the mean reverting constant is $\beta=10$ and the volatility constant is $\gamma=0.4$.

We can observe the same behavior that the one described in the previous simulations. However, the realizations of the solution of the random equation \eqref{logistical} are much closer to the deterministic ones since $\beta$ is larger.\n

%\alain{
\section{Parameter estimation in the logistic model.}\label{s4}

We consider again the logistic equation that we write
\begin{eqnarray}
\frac{dx}{dt} & = & ax\left(1-\frac{x}{K}\right),
\end{eqnarray}
where we put $K=ac$. In several ecological systems, $K$ represents a
{\em carrying capacity}, which is related to the size of the
population when all the sites are colonized. When considering the density of the
population or the proportion $p=x/K \in [0,1]$ of occupied sites, the
variable $p$ is solution of the differential equation
\begin{eqnarray}
\label{dyn_p}
\frac{dp}{dt} & = & rp(1-p)
\end{eqnarray}
where the parameter $r=aK$, usually known as the {\em intrinsic growth
  rate}, may fluctuate about a nominal value under environmental
variations (season, light, temperature...).
We consider then random perturbations on $r$, as in previous section
\begin{eqnarray}
\frac{dp}{dt} & = & (r+\alpha z^*_{\beta,\gamma}(\theta_t\omega))p(1-p) ,\label{pr}
\end{eqnarray}
for which there exist positive numbers $\underline{r}$, $\bar r$ with
$\underline{r}<\bar r$, such that each realization of $r+\alpha
z^*_{\beta,\gamma}(\theta_t\omega)$ belongs to the interval $(\underline
r,\bar r)$.
The question under investigation in this section is the estimation of
the parameter $r$, measuring the proportion $p$ over the time, in both
deterministic and random frameworks. 

\medskip 

\begin{figure}[H]
\includegraphics[scale=0.25]{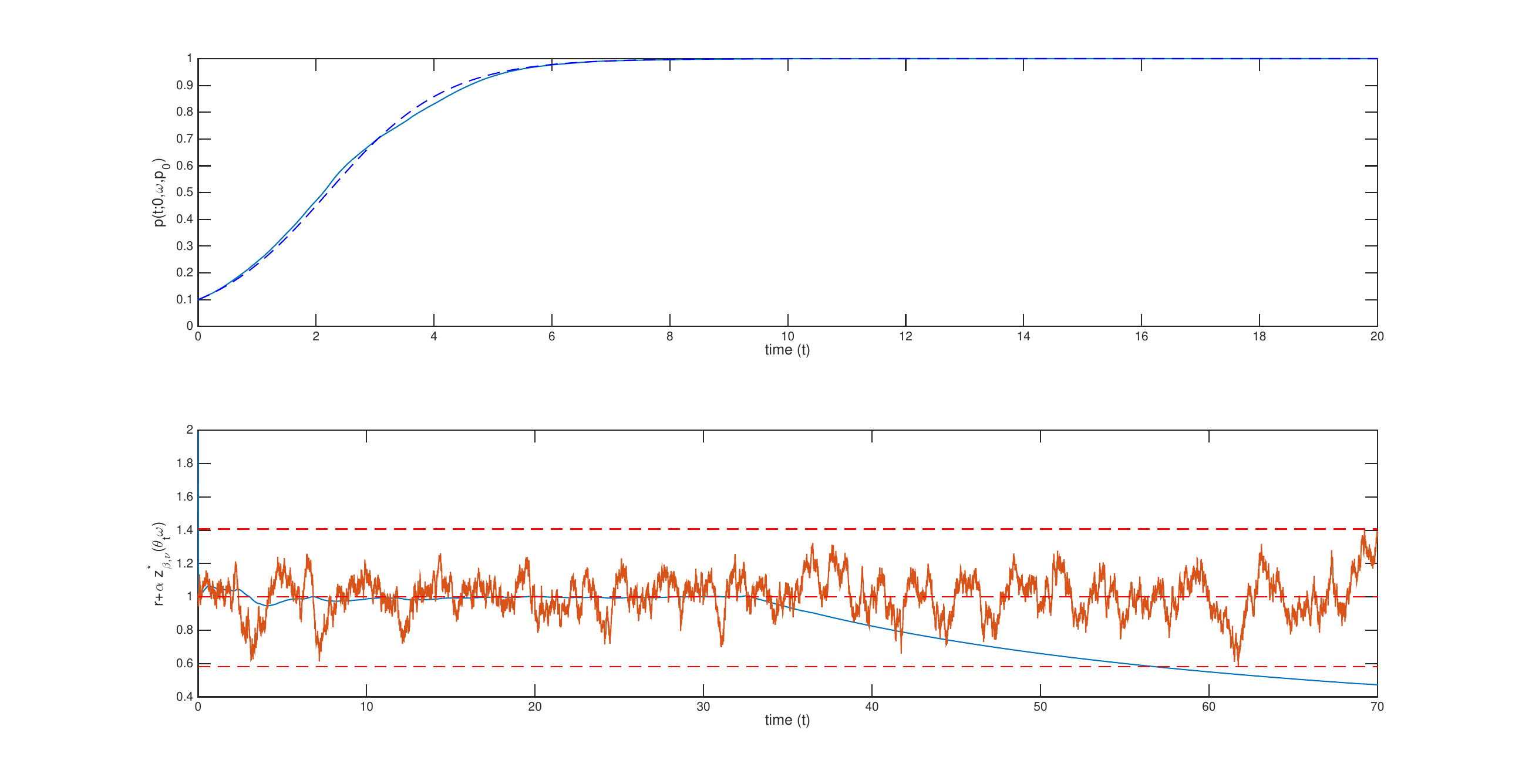}
\caption{\label{figsimurestim} Realizations of the solution of \eqref{pr} and estimation of $r$ in presence of random perturbation.}
\end{figure}

From equation \eqref{dyn_p}, one obtains an exact
expression of $r$, when the dynamics of $p$ is not at steady state
\begin{equation}
r=\frac{1}{t}\int_{p(0)}^{p(t)} \frac{dp}{p(1-p)}
=\frac{1}{t}\log\left(\frac{p(t)(1-p(0))}{p(0)(1-p(t))}\right), \quad t>0,
\end{equation}
One can then consider
\begin{equation}
\label{restim}
\hat r(t)=\frac{1}{t}\log\left(\frac{p(t)(1-p(0))}{p(0)(1-p(t))}\right), \quad t>0
\end{equation}
as an estimator of $r$ in presence of random
perturbation. Simulations show that this estimator behaves well as long
as $p(t)$ is not too close from its limiting value $1$ (see Figure \ref{figsimurestim}, where two panels are presented: the first shows the realizations of the solution of \eqref{pr} and the other overlaps the realization $r+\alpha z^*_{\beta,\gamma}(\theta_t\omega)$ (orange line), its estimator (blue line) and the dashed lines represent the values $\underline{r}$, $r$ and $\bar{r}$.).
Differently to a true observer \alain{(see for instance \cite{GK01} for an introduction to the theory of observers)}, we have no information to know when to
trust this estimator. 

\medskip

\begin{figure}[H]
\includegraphics[scale=0.25]{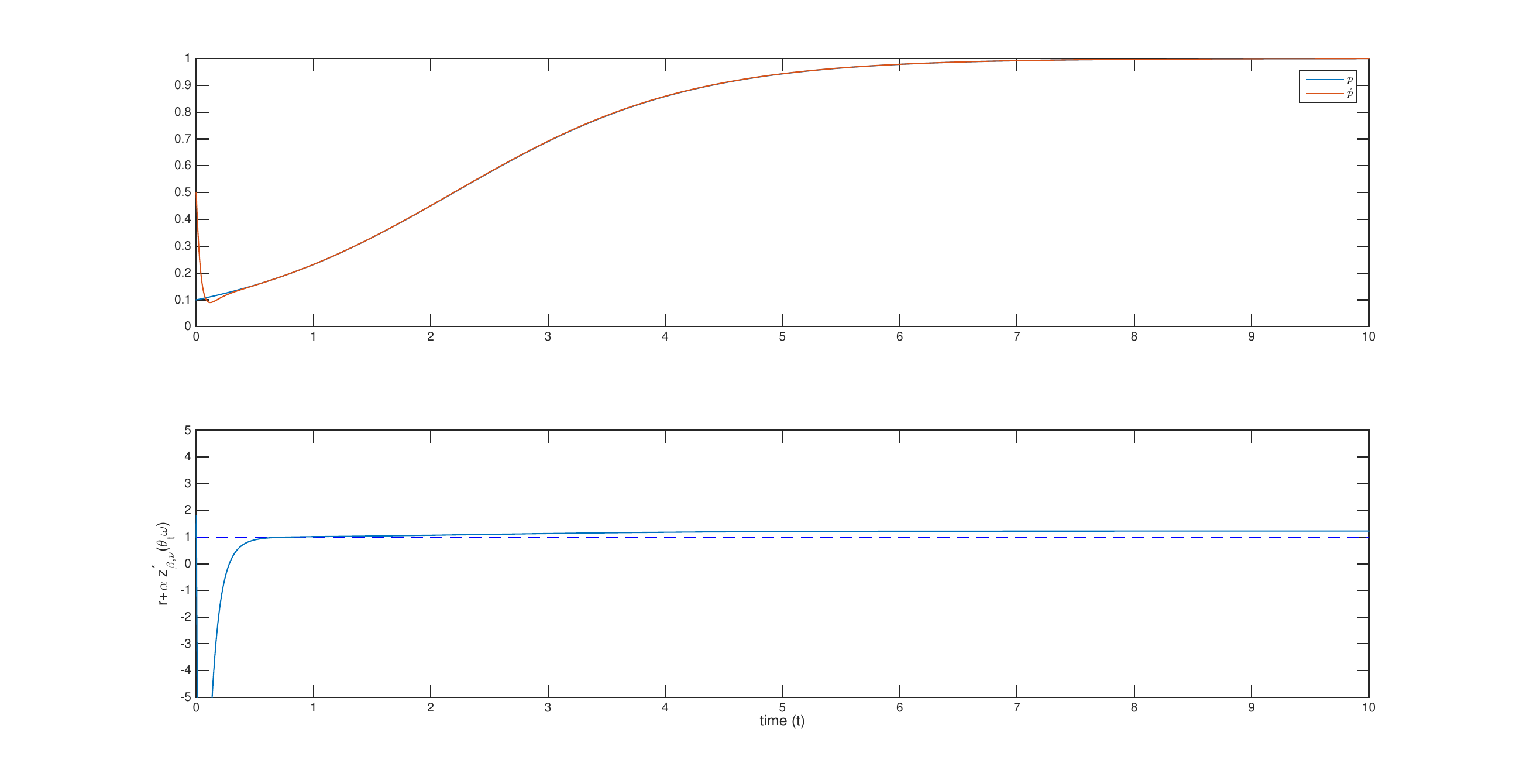}
\caption{\label{figsimuNormalObs} Simulation of the observer in normal
form for $T=2$ and $\theta=15$}
\end{figure}

We look instead for observers, i.e.~dynamical systems with on-line corrective terms (see Appendix \ref{appendix} for the main ingredients which are not here to avoid technicalities). After setting up an observer in normal form, it is possible to obtain a practical
convergence of the observer, in the sense that for any
$\varepsilon>0$, \alain{one can design an observer such that there exits $T>0$ with} $\hat
r(T) \in [r-\varepsilon,r+\varepsilon]$. One cannot expect a better convergence
of this observer as the system is not observable at $p=1$ and all 
solutions with non null $p(0)$ converge asymptotically to this
singular point.
On the simulations depicted on Figure \ref{figsimuNormalObs}, one can see two different panels: the first one shows the solution of \eqref{pr} without the presence of noise and the observer $\hat{p}$ and the second one represents the estimator $\hat{r}$ and $r$. Moreover, we can observer that the error of the observer
 has converged much before the system has reached the neighborhood of the
 steady state $p=1$\alain{, as it is desired for a true observer. Moreover the innovation, that is the difference between the observed variable $x$ and the variable $\hat x$ (in blue) of the observer informs on the convergence of the estimator $\hat{r}$ (when $\hat{x}$ stays almost equal to $x$, we know also that $\hat{r}$ stays close to the unknown value $r$).}
 
\begin{figure}[H]
\includegraphics[scale=0.25]{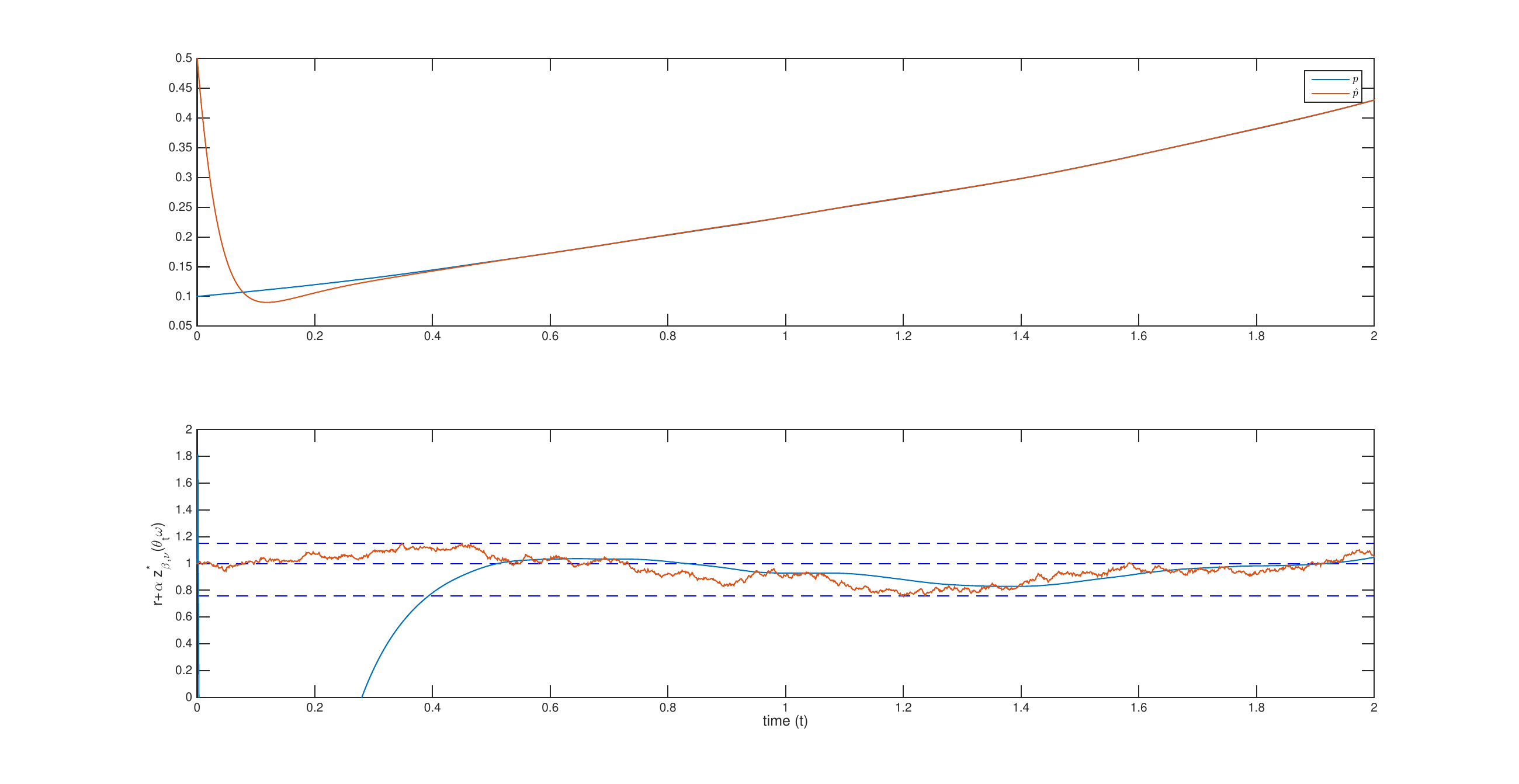}
\caption{\label{figsimuNormalObs+noise} Simulation of the observer in normal
form for $T=2$ and $\theta=15$ with random perturbation on $r$}
\end{figure}

Finally, in Figure \ref{figsimuNormalObs+noise} we present simulations of this observer in presence of random
perturbations. The first panel shows the solution of \eqref{pr} in presence of noise by means of the O-U process and the second panel overlaps the realization of $r+\alpha z^*_{\beta,\gamma}(\theta_t\omega)$ (orange line) and the estimator $\hat{r}$ (blue line). The dashed lines represent the values of $\bar{r}$, $r$ and $\underline{r}$.
%}

\section{Random competitive Lotka-Volterra models}\label{s5}

In this section we consider a competitive Lotka-Volterra model given by
\begin{eqnarray}
\f{dx}{dt}&=&x(\lambda-a x-by),\label{d1}
\\[1.3ex]
\f{dy}{dt}&=&y(\mu -cx-ey)\label{d2}
\end{eqnarray}
\noindent where $x=x(t)$ is the number of population of the {\it first} species, $y=y(t)$ is the number of population of the {\it second} species, $\lambda$ and $\mu$ are the specific growth rates of each species, respectively, $a$ and $c$ are the carrying capacities of each species, respectively, $b$ measures the interaction of the {\it first} species on the {\it second} one and $e$ is the interaction that the {\it second} species has on the {\it first} ones. We assume that all the parameters in the system are positive.\n

Concerning the deterministic competitive system \eqref{d1}-\eqref{d2}, it is known that coexistence of both populations can be ensured as long as conditions
\begin{equation}
\f{d}{b}>\f{\mu}{\lambda}>\f{c}{a}\quad \text{and}\quad ad-bc>0
\end{equation}
\noindent hold true.\n

In this case we are interested in studying the previous system where the growth rates are affected by the O-U process. Then, we consider the random competitive model given by
\begin{eqnarray}
\f{dx}{dt}&=&x(\lambda+\alpha z^*_{\beta,\gamma}(\theta_t\omega)-a x-by),\label{r1}
\\[1.3ex]
\f{dy}{dt}&=&y(\mu+\alpha z^*_{\beta,\gamma}(\theta_t\omega) -cx-ey)\label{r2}
\end{eqnarray}
\noindent where $z^*_{\beta,\gamma}(\theta_t\omega)$ denotes the O-U process introduced in Section \ref{secou} of this work and $\alpha>0$ represents the intensity of the noise.\n

Thanks to a suitable choice of the parameter $\beta$, we know that $\lambda+\alpha z^*_{\beta,\gamma}(\theta_t\omega)\in [\underline{\lambda}, \bar{\lambda}]$ and $\mu+\alpha z^*_{\beta,\gamma}(\theta_t\omega)\in [\underline{\mu}, \bar{\mu}]$ for any $t\in\R$.\n

Hence, from the random competitive system \eqref{r1}-\eqref{r2}, we can obtain the following differential inequalities for the dynamics of the population of both species involved in our model
\begin{eqnarray}
\f{dx}{dt}&\leq&x(\bar{\lambda}-ax),\label{a1}
\\[1.3ex]
\f{dy}{dt}&\leq&y(\bar{\mu}-ey).\label{a2}
\end{eqnarray}

Then, we obtain that the population of the both species are bounded from above
\begin{equation}
x(t;0,\omega,x_0)\leq\f{\bar{\lambda}}{a}\quad \text{and}\quad y(t;0,\omega,y_0)\leq \f{\bar{\mu}}{e}
\end{equation}
\noindent for any initial values $x_0\geq 0$, $y_0\geq 0$, any $\omega\in\Omega$ and $t\geq 0$.\n

In addition, from \eqref{a1} it is possible to obtain the following differential inequality
\begin{equation}
\f{dx}{dt}\geq x\left(\underline{\lambda}-b\f{\bar{\mu}}{e}-ax\right)
\end{equation}
\noindent whence we can obtain its explicit solution which is given by
\begin{equation}
x(t;0,\omega,x_0)\geq\f{x_0}{e^{-\left(\underline{\lambda}-b\f{\bar{\mu}}{a}\right)t}+\f{x_0a}{\underline{\lambda}-b\f{\bar{\mu}}{e}}\left(1-e^{-\left(\underline{\lambda}-b\f{\bar{\mu}}{a}\right)t}\right)}
\end{equation}
\noindent for any initial value $x_0\geq 0$, any $\omega\in\Omega$ and $t\geq 0$, from which, by taking limit when $t$ goes to infinity, we obtain
\begin{equation}
\lim_{t\rightarrow+\infty}x(t;0,\omega,x_0)\geq\f{\underline{\lambda}-b\f{\bar{\mu}}{e}}{a}.
\end{equation}

Then, the population of the {\it first} species persists as long as the condition
\begin{equation}
\f{\bar{\mu}}{\underline{\lambda}}<\f{e}{b}\label{e1}
\end{equation}
\noindent is fulfilled.\n

Concerning the other population, the same argument can be done and we obtain
\begin{equation}
\lim_{t\rightarrow+\infty}y(t;0,\omega,y_0)\geq\f{\underline{\mu}-\f{\bar{\lambda}}{a}}{e},
\end{equation}
\noindent then the population of the {\it second} species persists as long as the following condition is fulfilled
\begin{equation}
\f{\underline{\mu}}{\bar{\lambda}}>\f{c}{a}.\label{e2}
\end{equation}

In conclusion, for any $\varepsilon>0$, $\omega\in\Omega$ and every initial values $x_0\geq 0$ and $y_0\geq 0$, there exists some time $T(\varepsilon,\omega)>0$ such that the solution of the random system \eqref{a1}-\eqref{a2} can be bounded inside the frame
\begin{equation}
B_\varepsilon=\left[\f{\underline{\lambda}-b\f{\bar{\mu}}{e}}{a}-\varepsilon,\f{\bar{\lambda}}{a}\right]\times\left[\f{\underline{\mu}-\f{\bar{\lambda}}{a}}{e}-\varepsilon,\f{\bar{\mu}}{e}\right]
\end{equation}
\noindent for every $t\geq T(\varepsilon,\omega)$.\n

Therefore, for any $\varepsilon>0$, $B_\varepsilon$ is a strictly positive deterministic absorbing set for the solutions of the system \eqref{a1}-\eqref{a1}, whence we have that
\begin{equation}
B_0=\left[\f{\underline{\lambda}-b\f{\bar{\mu}}{e}}{a},\f{\bar{\lambda}}{a}\right]\times\left[\f{\underline{\mu}-\f{\bar{\lambda}}{a}}{e},\f{\bar{\mu}}{e}\right]\label{ab}
\end{equation}
\noindent is a strictly positive deterministic attracting set for the solutions of the system \eqref{a1}-\eqref{a2}.\n

From the previous analysis, we can observe that, as long as conditions \eqref{e1} and \eqref{e2} are satisfied, we can ensure the coexistence of the population of both species.\n

\begin{figure}[H]
\begin{center}
\includegraphics[scale=0.25]{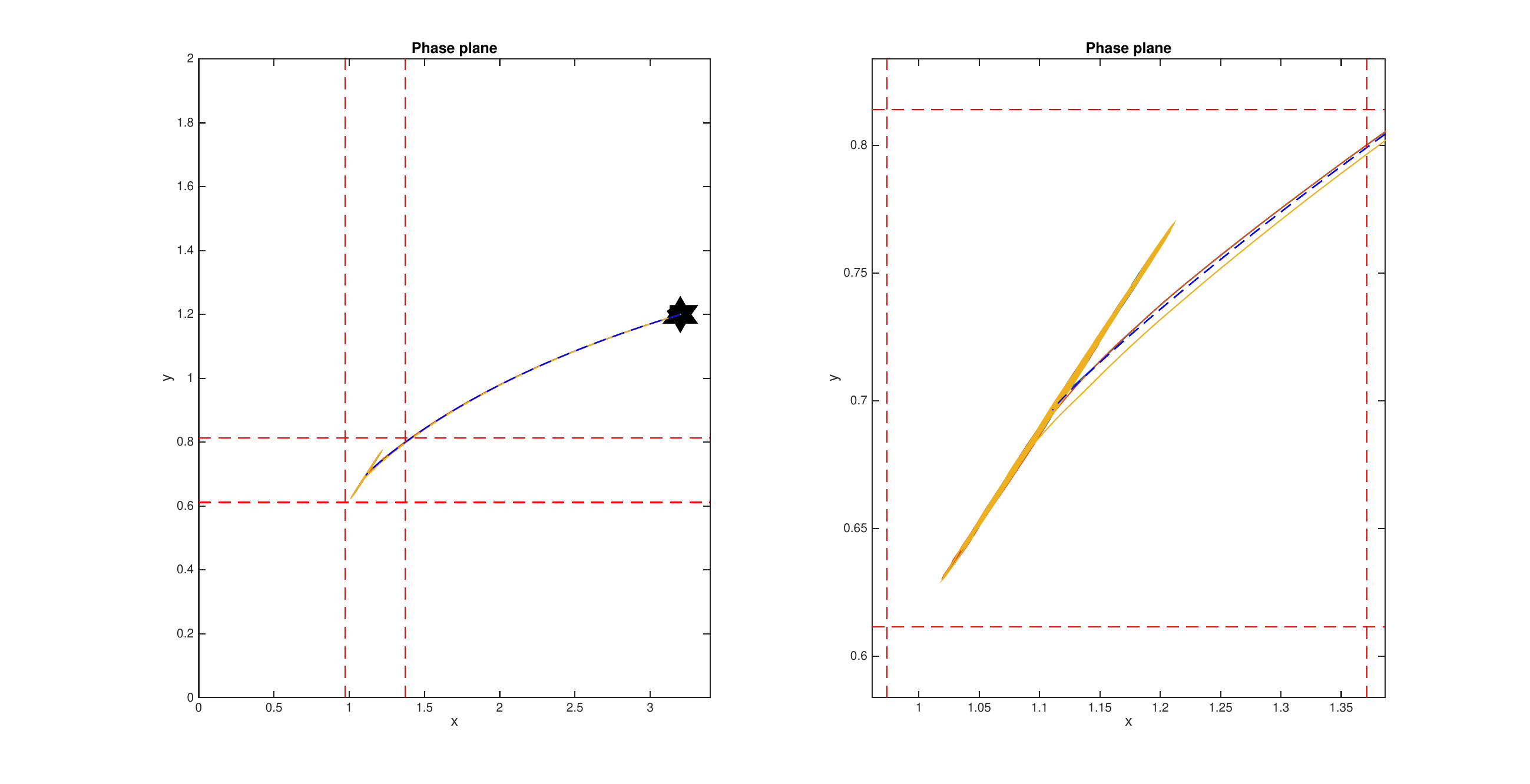}
\end{center}
\caption{Phase plane with realizations of the solutions of the competitive Lotka-Volterra system for $x_0=3.2$ and $y_0=1.2$}
\label{c0}
\end{figure}

Now, we present some numerical simulations to illustrate the results provided in this section. In Figure \ref{c0} we show the phase plane with several realizations of the solution of the random competitive system \eqref{a1}-\eqref{a2} for the initial values $x_0=3.2$ and $y_0=1.2$ and the following values of the rest of the parameters $a=20$, $b=4$, $c=1$, $e=30$, $\lambda=25$, $\mu=22$, the amount of noise is $\alpha=2$, the mean reverting constant is $\beta=1$ and the volatility constant is $\gamma=0.5$. We remark that the right panel shows a zoom of the left one to see the absorbing set of the solutions, which is the box delimited by the dashed lines.

We can observe that all the realizations of the solution of the system remain, after some time, inside a rectangle limited by the red dashed lines. This rectangle is the absorbing set $B_0$ (see \eqref{ab}) obtained in the mathematical results which is deterministic.\n

\begin{figure}[H]
\begin{center}
\includegraphics[scale=0.25]{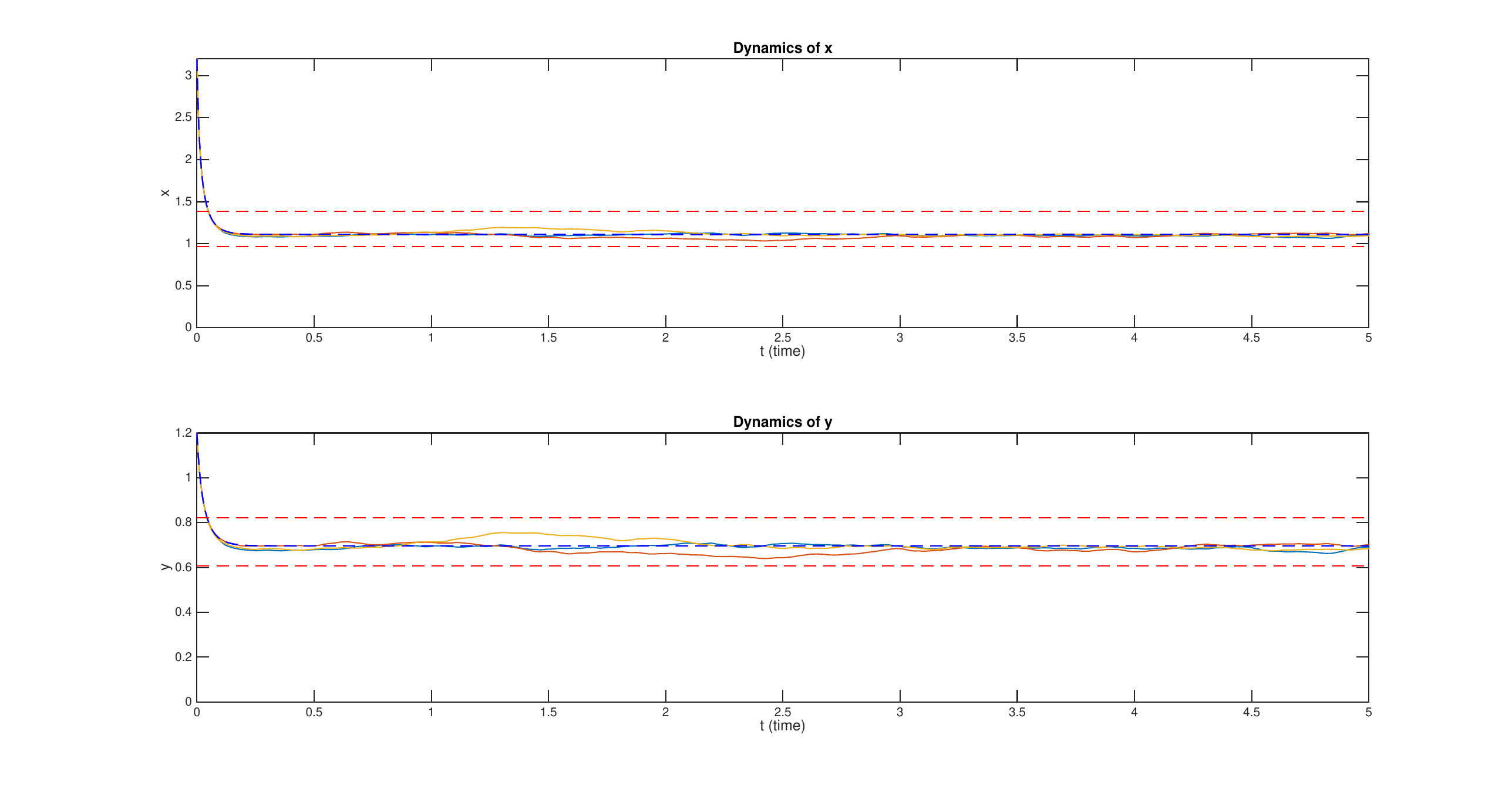}
\end{center}
\caption{Realizations of the solution of the competitive Lotka-Volterra system (both state variables depending on time) for for $x_0=3.2$ and $y_0=1.2$}
\label{c02}
\end{figure}

In Figure \ref{c02} we present the dynamics of both species individually where the red dashed lines represent the bounds guaranteed for the corresponding state variables. We can observe that both species are fluctuating around the deterministic solution inside a strictly positive interval that allows us to guarantee the persistence of both species. In addition, these intervals are deterministic \alain{in the sense that} they do not depend on the realization on the noise and can be chosen as explained previously.

\begin{figure}[H]
\begin{center}
\includegraphics[scale=0.25]{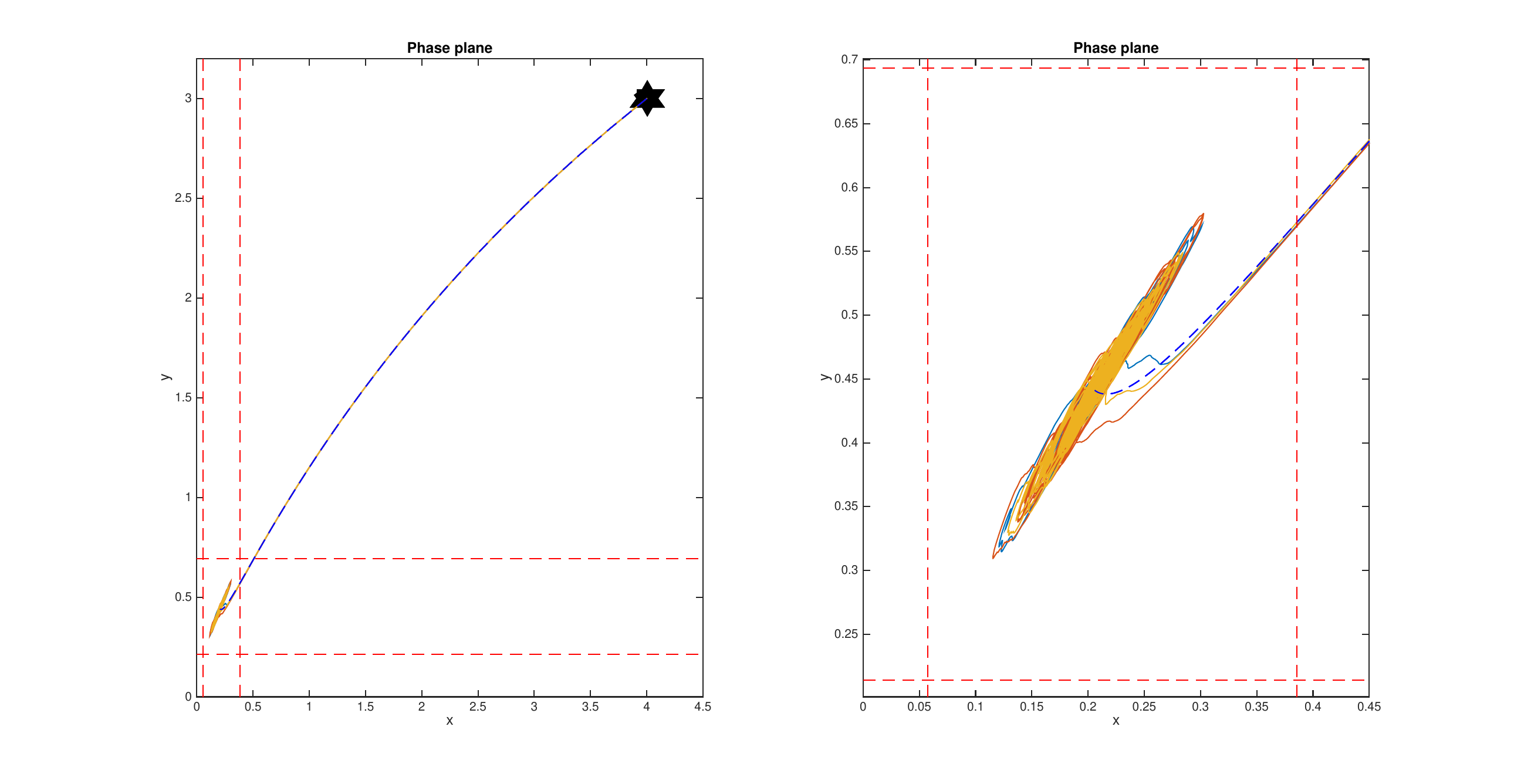}
\end{center}
\caption{Phase plane with realizations of the solutions of the competitive Lotka-Volterra system for $x_0=4$ and $y_0=3$}
\label{c2}
\end{figure}

In Figure \ref{c2} we present the phase plane with several realizations of the solutions of the random competitive system \eqref{a1}-\eqref{a2} for the initial values $x_0=4$ and $y_0=3$ and the following values of the rest of the parameters $a=20$, $b=2$, $c=4$, $e=314$, $\lambda=5$, $\mu=7$, the amount of noise is $\alpha=2$, the mean reverting constant is $\beta=1$ and the volatility constant is $\gamma=0.5$. \alain{Let us underline} that the right panel is a zoom of the left one to show the absorbing set $B_0$ (see \eqref{ab}), which is the box bounded by the dashed lines.

\begin{figure}[H]
\begin{center}
\includegraphics[scale=0.25]{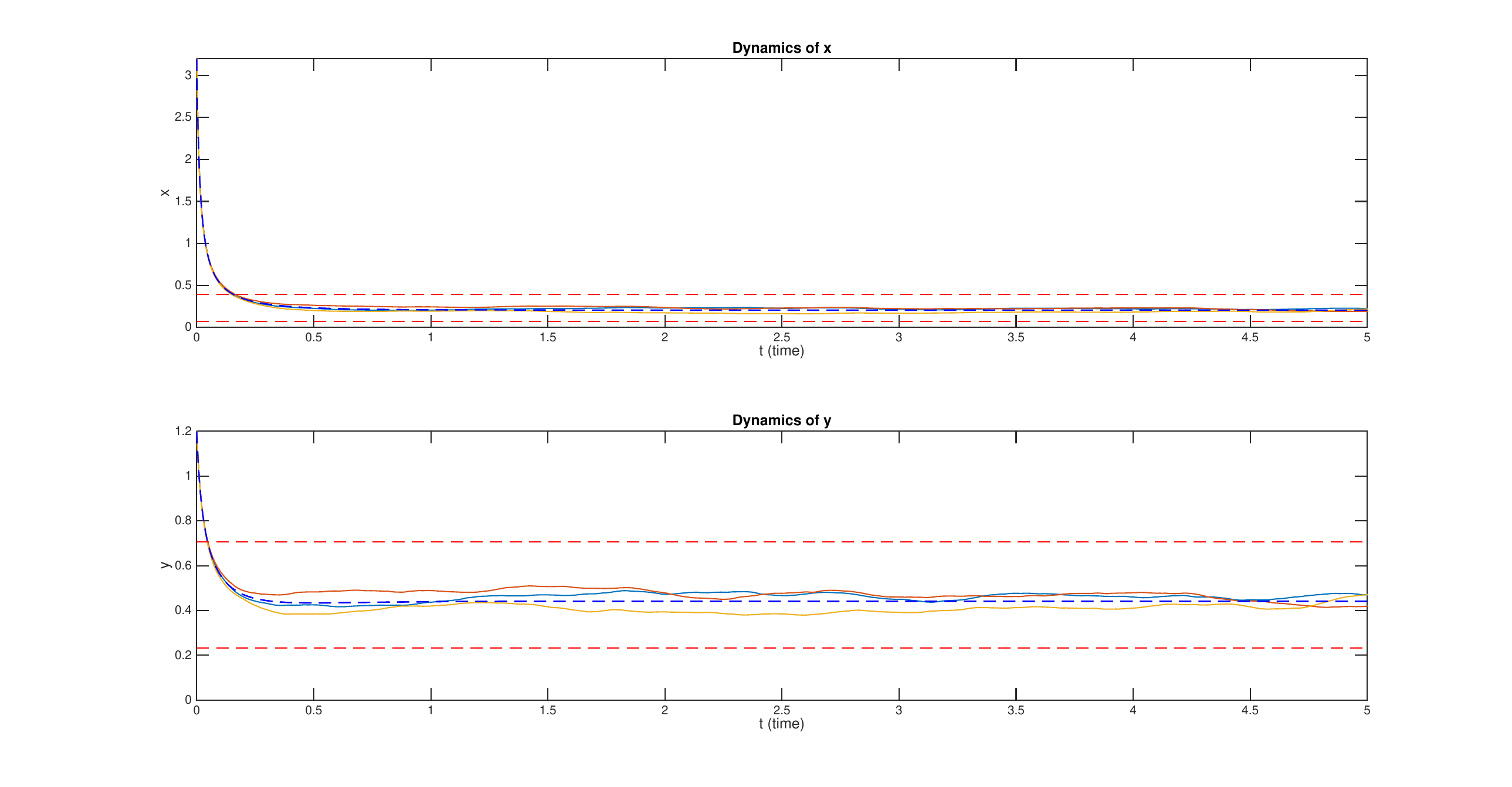}
\end{center}
\caption{Realizations of the solution of the competitive Lotka-Volterra system (both state variables depending on time) for $x_0=4$ and $y_0=3$}
\label{c22}
\end{figure}

In Figure \ref{c22} we can observe the dynamics of both species individually. In this case we have changed the parameters of the system and we can see that all the realizations of the solution of the random system remain inside the strictly positive absorbing set delimited by the red dashed lines, as proved previously, which means that we can ensure the persistence of both species also in this case.

\section{Random chemostat model}\label{s6}

In this short section we would like to \alain{discuss about} the O-U process \alain{that has already being considered for modelling random perturbations on the input flow rate} in the chemostat model. Every detail about the way of modeling and a complete analysis of the resulting random model can be found in \cite{CGLR,JLC} thus we will omit the details in this section. Here, we just give some comments concerning the work in \cite{CGLR,JLC}.\n

As already explained in the introductory section, several drawbacks can be found when perturbing the input flow of the chemostat model by using the standard Wiener process (see \cite{CGL,corrigendumchapter}). For instance, the input flow could take extremely large values \alain{and thus could } negative. Due to this fact, which is unrealistic from the biological point of view since we know that the input flow is fluctuating in a positive bounded interval, we also have that some state variables \alain{that describe population size could take negative values} which is also unrealistic from the biological point of view. In addition, it is not posible to ensure the persistence of the species \alain{as we did and which corresponds to real observations}.\n

However, everything \alain{these drawbacks are circumvented} when introducing the perturbations on the input flow by means of the O-U process as explained in this paper. The first important improvement is the fact that the perturbed input flow is ensured to be bounded, as in \alain{real experiments} (see Figure \ref{real} where we presented the real data). In addition, it is possible to prove that there exists absorbing and attracting sets which are deterministic (then they do not depend on the realization of the noise) and moreover, what is essential to prove the persistence of the species, positive. Furthermore, these results are proved in forward sense which \alain{suits the point of view of applications.}.

\section{Conclusions and final comments}\label{s7}

In this final section we would like to \alain{draw} some conclusions and final comments concerning the O-U process introduced in Section \ref{secou}. We recall that the most important improvement of \alain{this way of modeling the noise, compared to} other kinds of noise considered previously in the literature is the fact that the O-U process depends of two parameters, the volatility constant $\gamma$ and the mean reversion constant $\beta$ which play an important role and allow the noise to have all the expected properties that we have formulated in the introduction.\n

In order to show the \alain{relevance} of this new \alain{way of modeling } we have presented in the previous sections some examples which illustrate the effect of this bounded noise when perturbing some very well-known models such as the logistic \alain{growth} or the Lotka-Volterra \alain{competition}. In addition, in \cite{CGLR,JLC} the authors consider this noise to model random input flows in the chemostat model where some relevant improvements are also achieved. Finally, this way of modeling is full of advantages \alain{from the mathematical analysis point of view but also,} which is essential, \alain{as a quite realistic modeling} from the biological point of view.\n

In conclusion, \alain{we believe that this modeling approach is generic and could be applied in most of the population models, when some aspects or parameters are expected to subject to randomness with bounded realizations}. For instance, it could be very interesting to analyze prey-predator models where, in the deterministic case periodic orbits and limit cycles are present. In this way, it could be possible to define a concept of {\it random} periodic orbits in the sense that the solutions of the system are fluctuating around the deterministic periodic orbit inside some interval that depends on the parameters of the O-U process, as in the examples of the present paper. Another idea is to analyze the problem of the observer with measurements perturbed by the O-U process. These are some ideas among other ones to carry on  the applications of this way of modeling noise.

\appendix

\section*{Appendix}\label{appendix}

\begin{figure}[H]
\includegraphics[scale=0.25]{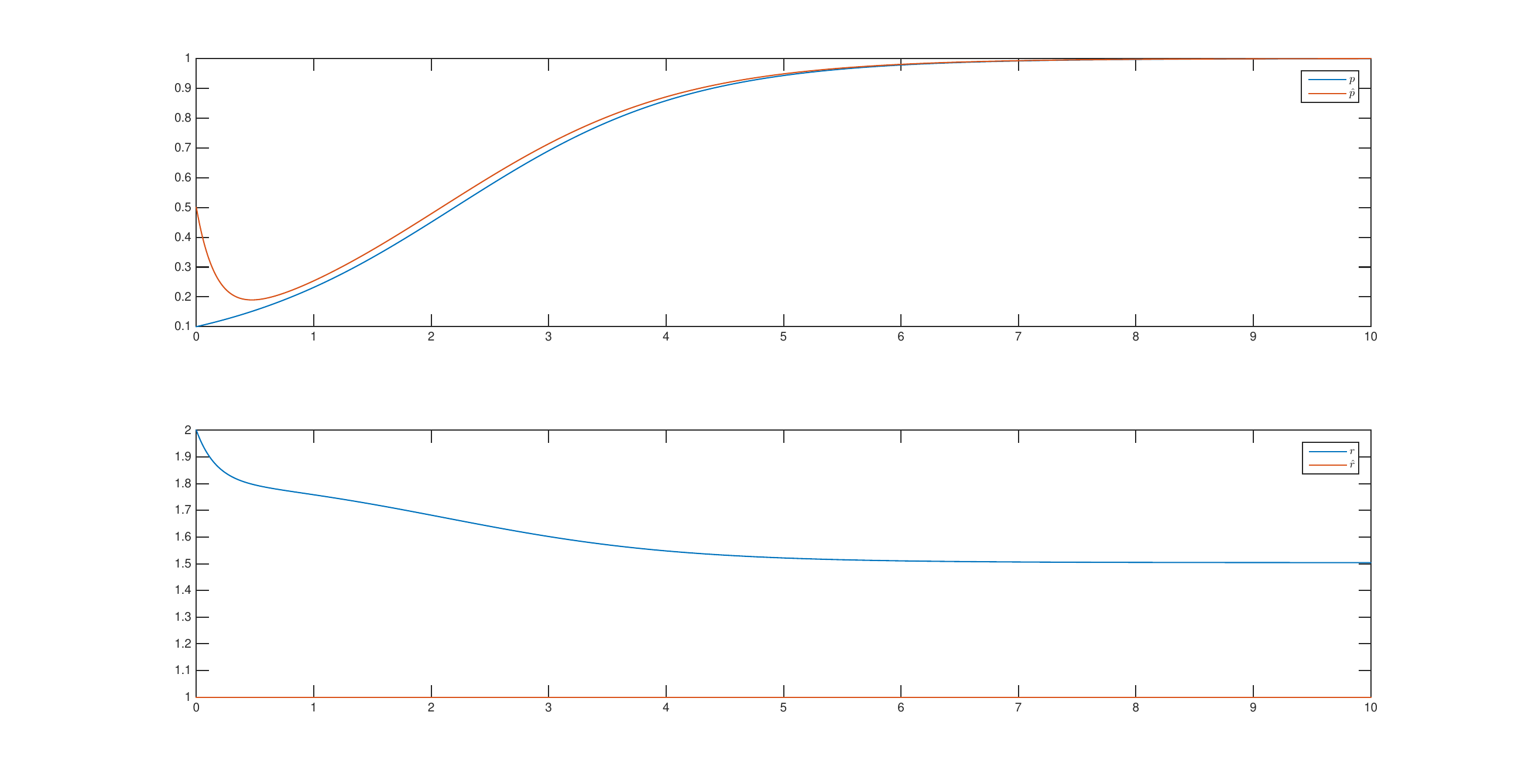}
\caption{\label{figsimuLuenberger} Simulation of the Luenberger
  observer with $\gamma_{}=-5$ and $\gamma_{r}=-1$}
\end{figure}

\alain{Let us give more details about the observer construction used in Section \ref{s4}.}
\alain{For simplicity, let us first consider} by the deterministic framework. For this
purpose, we consider the extended dynamics
\begin{eqnarray}
\frac{dp}{dt} & = & rp(1-p) \label{dp} \\
\frac{dr}{dt} & = & 0 \label{dr}
\end{eqnarray}
with the measured output
\begin{equation*}
y(y)=p(t)
\end{equation*}
Notice this system is not observable \alain{(see \cite{GK01} for the definition of observability)} at the steady states $p=0$ or
$p=1$. 
When the system is not at equilibrium, let us first consider a classical observer of Luenberger form
\begin{eqnarray}
\nonumber
\f{d\hat p}{dt}&=& \hat r y(t)(1-y(t))+G_{1}(\hat p -y(t))
\\[1.3ex]
\nonumber
\f{d\hat r}{dt}&=&G_{2}(\hat p-y(t))
\end{eqnarray}
where the gains parameters $G_{1}$, $G_{2}$ have to be chosen.
The dynamics of the error variables $e_{p}=\hat p-p$, $e_{r}=\hat r-r$
are given by the linear non-autonomous system
\begin{eqnarray}
\nonumber
\f{de_{p}}{dt}&=& G_{1}e_{p}+y(t)(1-y(t))e_{r}
\\[1.3ex]
\nonumber
\f{de_r}{dt}&=&G_{2}e_{p}
\end{eqnarray}
Consider then the quadratic function
\begin{equation}
\label{fctV}
V(e_{p},e_{r}) = \frac{1}{2}(e_{p}+\gamma_{r} e_{r})^2 + \frac{1}{2} e_{r}^2
\end{equation}
where $\gamma_{r}$ is a parameter. Notice that $V$ is definite positive
for any value of $\gamma_{r}$. One has, along any trajectory:
\begin{eqnarray}
\nonumber
\frac{dV}{dt} & = &
(e_{p}+\gamma_{r}e_{r})(G_{1}e_{p}+y(t)(1-y(t))e_{r}+\gamma_{r} G_{2}e_{p})
+e_{r}G_{2}e_{p}\\
\nonumber
 & = & (G_{r}+\gamma_{r} G_{2})e_{p}^2 + \gamma_{r} y(t)(1-y(t))e_{r}^2 +
 (\gamma_{r} G_{1}+\gamma_{r}^2G_{2}+G_{2})e_{p}e_{r}
\end{eqnarray}
Take $G_{1}<0$ and $\gamma_{r}<0$ and set
\begin{equation*}
G_{2}=-\frac{\gamma_{r}}{1+\gamma_{r}^2}G_{1}
\end{equation*}
Notice that for such choice, one has
\begin{equation*}
\gamma_{p}:=G_{1}+\gamma_{r} G_{2}=\frac{G_{1}}{1+\gamma_{r}^2}<0
\end{equation*}
Equivalently, $G_{1}$ and $G_{2}$ are defined as
\begin{eqnarray}
\nonumber
G_{1} & = & (1+\gamma_{r}^2)\gamma_{p}\\
\nonumber
G_{2} & = & -\gamma_{r} \gamma_{p}
\end{eqnarray}
with $\gamma_{p}$ and $\gamma_{r}$ negative.
For $V>0$, one has the inequality
\begin{equation*}
\frac{dV}{dt} = \gamma_{p} e_{p}^2 + \gamma_{r}
y(t)(1-y(t))e_{r}^2 < 0 , \quad \forall t>0.
\end{equation*}
However, we cannot conclude about the convergence of $V$ to $0$ because
\[
\int_{0}^{+\infty} y(t)(1-y(t))dt =\frac{1-p(0)}{r}< + \infty
\]
as this is shown on Figure \ref{figsimuLuenberger}.

\medskip

Consider now a second kind of observer, but in {\em normal form} \alain{(see \cite{GHO92,GK01})}, which
consists in applying the change of coordinates $(p,r) \to (z_{1},z_{2})$ with
\begin{eqnarray}
\nonumber
z_{1} & = & p\\
\nonumber
z_{2} & = & rp(1-p)
\end{eqnarray}
When $p$ is different from the steady states $p=0$ and $p=1$,
parameter $r$ can be reconstructed as
\begin{equation*}
r=\varphi(z_{1},z_{2}):=\frac{z_{2}}{z_{1}(1-z_{1})}
\end{equation*}
Dynamics \eqref{dp}-\eqref{dr} in these coordinates writes as follows
\begin{equation*}
\f{dz}{dt} = \underbrace{\left[\begin{array}{cc}0 & 1\\ 0 & 0
\end{array}\right]}_{A}z+\left[\begin{array}{c}0\\ \psi(z_{1},z_{2})
\end{array}\right] \quad \mbox{with} \quad
\psi(z_{1},z_{2}):=\varphi(z_{1},z_{2})\left(1-\frac{z_{1}}{2}\right)z_{2}
\end{equation*}
with the observation
\begin{eqnarray}
\nonumber
y(t) & = & \underbrace{\left[\begin{array}{cc} 0 & 1\end{array}\right]}_{C}z(t)
\end{eqnarray}
This leads to consider the following observer
\begin{eqnarray}
\nonumber
\f{d\hat z_{1}}{dt}&=& \hat z_{2}+G_{1}(\hat z_{1}-y(t))\\[1.3ex]
\nonumber
\f{d\hat z_{2}}{dt}&=& \psi(y(t),\hat z_{2}) + G_{2}(\hat z_{1}-y(t))
\end{eqnarray}
with the estimator
\begin{equation*}
\hat r(t) = \varphi(y(t),\hat z_{2}(t))
\end{equation*}
Notice that when the system is not at steady state, $\varphi(y(t),\hat
z_{2}(t))$ and $\psi(y(t),\hat z_{2}(t))$ are well defined for any
$t>0$.
The dynamics of the error $e=\hat z-z$ is given by the system
\begin{equation*}
\f{de}{dt}=\left[\begin{array}{cc}G_{1} & 1\\ G_{2} & 0
\end{array}\right]e+\left[\begin{array}{c}0\\ 1
\end{array}\right](\psi(y(t)-\hat z_{2})-\psi(y(t)-z_{2}))
\end{equation*}
The map $z_{2} \mapsto \psi(y(t),z_{2})$ is not Lipschitz
with respect to $z_{2}$ uniformly w.r.t.~$t$. However for any fixed $T$, it is
Lipschitz on any compact set uniformly on $[0,T]$. We can then use the
theory of {\em high-gains observers} \cite{GHO92,GK01}, which guarantees an exponential
decrease of the norm of the error on $[0,T]$, when the gains $G_{1}$, $G_{2}$ are
chosen such that
\begin{equation*}
\left[\begin{array}{c}
G_{1}\\
G_{2}
\end{array}\right] = -S_{\theta}^{-1}C^{\top}
\end{equation*}
where $S_{\theta}$ is the symmetric definite positive matrix solution of the Lyapunov equation
\begin{equation*}
A^{\top}S_{\theta}+S_{\theta}A-C^{\top}C+\theta S_{\theta}=0
\end{equation*}
and parameter $\theta>0$ is large enough. On can check that this gives
\[
\left[\begin{array}{c}
G_{1}\\
G_{2}
\end{array}\right] =\left[\begin{array}{c}
-2\theta\\
-\theta^2
\end{array}\right]
\]

\bibliographystyle{elsarticle-num.bst}

\bibliography{CCLR_submitted}

\end{document}